\documentclass[11pt]{article}
\usepackage{amsmath, upgreek}
\usepackage{amssymb}
\usepackage{color}
\usepackage{amscd}
\usepackage{xspace}
\usepackage{verbatim}

\newcommand{\R}{\mathbb{R}}

\newcommand{\N}{\mathbb{N}}

\newcommand{\beq}{\begin{equation} }
\newcommand{\eqq}{\end{equation} }
\newcommand{\cuad}{{\sqcap\kern-.68em\sqcup}}
\newcommand{\abs}[1]{\mid #1 \mid}
\newcommand{\norm}[1]{\|#1\|}

\newtheorem{definition}{Definition}[section]
\newtheorem{teo}{Theorem}[section]

\newtheorem{proposition}{Proposition}[section]

\newtheorem{lemma}{Lemma}[section]
\newtheorem{corollary}{Corollary}[section]
\newtheorem{remark}{Remark}[section]
\newcommand{\bremark}{\begin{remark} \em}
\newcommand{\eremark}{\end{remark} }

\def\beeq{\begin{equation}}
\def\eeq{\end{equation}}
\newcommand{\begeqaet}{\begin{eqnarray*}}
\newcommand{\eneqaet}{\end{eqnarray*}}

\newcommand{\myfrac}[2]{{\displaystyle \frac{#1}{#2} }}
\newcommand{\myint}[2]{{\displaystyle \int_{#1}^{#2}}}

\begin{document}

\begin{center}{\bf  \Large  Fractional heat equations with subcritical absorption \\[2mm]
having a measure as initial data}\medskip


\bigskip
{\small
       {\bf Huyuan Chen}\footnote{chenhuyuan@yeah.net}\smallskip

       Department of Mathematics, Jiangxi Normal University,
       \\ Nanchang 330022, China
       \\[2mm]
     {\bf Laurent V\'{e}ron}\footnote{Laurent.Veron@lmpt.univ-tours.fr}

\smallskip
Laboratoire de Math\'{e}matiques et Physique Th\'{e}orique
\\  Universit\'{e} Fran\c{c}ois Rabelais, Tours, France\\[2mm]
{\bf Ying Wang}\footnote{yingwang00@126.com}

\smallskip
Departamento de Ingenier\'{\i}a  Matem\'atica
\\ Universidad de Chile, Chile\\[1mm]
}

\begin{abstract}
We study  existence and uniqueness of   weak solutions to (F) $\partial_t u+ (-\Delta)^\alpha
u+h(t, u)=0 $ in $(0,\infty)\times\R^N$,
with initial condition $u(0,\cdot)=\nu$ in $\R^N$, where $N\ge2$, the operator $(-\Delta)^\alpha$
is the fractional Laplacian with $\alpha\in(0,1)$, $\nu$ is
a bounded Radon measure and $h:(0,\infty)\times\R\to\R$ is a continuous function   satisfying a subcritical integrability condition.

In particular, if  $h(t,u)=t^\beta u^p$ with $\beta>-1$ and $0<p<p^*_\beta:=1+\frac{2\alpha(1+\beta)}{N}$, we prove that there exists a unique weak solution $u_k$ to (F) with $\nu=k\delta_0$, where $\delta_0$ is the Dirac mass at the origin. We obtain that  $u_k\to\infty$ in $(0,\infty)\times\R^N$ as $k\to\infty$ for $p\in(0,1]$ and the limit of $u_k$ exists as $k\to\infty$ when  $1<p<p^*_\beta$,
 we   denote it by $u_\infty$.
When  $1+\frac{2\alpha(1+\beta)}{N+2\alpha}:=p^{**}_\beta<p<p^*_\beta$,
$u_\infty$ is the minimal self-similar solution of $(F)_\infty$ $\partial_t u+ (-\Delta)^\alpha u+t^\beta u^p=0 $ in $(0,\infty)\times\R^N$  with the initial condition $u(0,\cdot)=0$ in $\R^N\setminus\{0\}$ and it  satisfies $u_\infty(0,x)=0$ for $x\neq 0$.
While  if $1<p<p^{**}_\beta$, then $u_\infty\equiv U_p$, where $U_p$ is the maximal solution of the differential equation  $y'+t^\beta y^p=0$ on $\R_+$.
\end{abstract}
\end{center}

\tableofcontents \vspace{1mm}

\vspace{1mm}


\vspace{2mm}

\setcounter{equation}{0}
\section{Introduction}
Let  $h:(0,\infty)\times\R\to\R$ be a continuous function and
$Q_\infty=(0,\infty)\times\R^N$ with $N\ge2$. The first object of this paper is to consider existence and uniqueness of weak solutions  to  fractional heat equations
\begin{equation}\label{he 1.1}
 \arraycolsep=1pt
\begin{array}{lll}
\partial_t u + (-\Delta)^\alpha  u+h(t,u)=0\quad & {\rm in}\quad Q_\infty,\\[2mm]
 \phantom{\partial_t u + (-\Delta)^\alpha u + \ \ }
u(0,\cdot)=\nu\quad & {\rm in}\quad \R^N,
\end{array}
\end{equation}
where    $\nu$ belongs to the space $\mathfrak{M}^b(\R^N)$ of bounded Radon measures in $\R^N$ and $(-\Delta)^\alpha $ ($0<\alpha<1$) is the fractional Laplacian defined by
$$(-\Delta)^\alpha  u(t,x)=\lim_{\epsilon\to0^+} (-\Delta)_\epsilon^\alpha u(t,x),$$
where, for $\epsilon>0$,
$$
(-\Delta)_\epsilon^\alpha  u(t,x)=\int_{\R^N}\frac{ u(t,x)-u(t,z)}{|z-x|^{N+2\alpha}}\chi_\epsilon(|x-z|) dz
$$
and
$$\chi_\epsilon(r)=\left\{ \arraycolsep=1pt
\begin{array}{lll}
0\qquad & {\rm if}\quad r\in[0,\epsilon],\\[2mm]
1\qquad & {\rm if}\quad r>\epsilon.
\end{array}
\right.$$

In a pioneering work, Brezis and Friedman \cite{BF} have studied the semilinear heat equation with measure as initial data
\begin{equation}\label{eq 08-10-13 1a}
 \arraycolsep=1pt
\begin{array}{lll}
\partial_t u -\Delta  u+ u^p=0\quad & {\rm in}\quad Q_\infty,\\[2mm]
 \phantom{ \partial_t\Delta  -\Delta   }
u(0,\cdot)=k\delta_0\quad & {\rm in}\quad \R^N,
\end{array}
\end{equation}
where $k>0$ and $\delta_0$ is the Dirac mass at the origin. They proved that if $1<p<(N+2)/N$, then for every $k>0$
there exists a unique solution $u_k$ to (\ref{eq 08-10-13 1a}). When $p\geq (N+2)/N$, problem (\ref{eq 08-10-13 1a}) has no solution and even more, they proved that no nontrivial solution of the above equation vanishing on $\R^N\setminus\{0\}$ at $t=0$ exists.
 When $1<p<1+\frac{2}{N}$,  Brezis, Peletier and Terman used a dynamical system technique in \cite{BPT}  to prove the existence of a
{\it very singular solution} $u_s$ to
\begin{equation}\label{eq 08-10-13 2a}
\arraycolsep=1pt
\begin{array}{lll}
\partial_t u -\Delta  u+u^p=0\quad & {\rm in}\quad Q_\infty,
\end{array}
\end{equation}
vanishing at $t=0$ on $\R^N\setminus\{0\}$. This function $u_s$ is self-similar, i.e. expressed under the form
\begin{equation}\label{self-Sym}
u_s(t,x)=t^{-\frac{1}{p-1}}f\left(\frac{|x|}{\sqrt t}\right),
\end{equation}
and $f$ is uniquely determined by the following conditions
\begin{equation}\begin{array}{lll}\label{zero-0}
f''+\left(\frac{N-1}{\eta}+\frac{1}{2}\eta\right)f'+\frac{1}{p-1}f-f^p=0\quad\text{on }\;\R_+\\[2mm]
f>0\quad\text{and }\ f\text{ is smooth on }\R_+\\[2mm]
f'(0)=0\quad\text{and }\;\lim_{\eta\to\infty}\eta^{\frac{2}{p-1}}f(\eta)=0.
\end{array}\end{equation}
Furthermore, it satisfies
$$ f(\eta)=c_1e^{-\eta^2}\eta^{\frac2{p-1}-N}\{1-O(|x|^{-2})\}\quad\text{as }\eta\to\infty $$
for some  $c_1>0$.
Later on,  Kamin  and Peletier in \cite{KP} proved that the sequence of weak solutions $u_k$
converges to the very singular solution $u_s$ as $k\to\infty$. After that, Marcus and V\'{e}ron in \cite{MV1}  studied the equation in the framework of the {\it initial trace} theory. They pointed out the role of the very singular solution of (\ref{eq 08-10-13 2a}) in the study of the singular set of the initial trace, showing in particular that it is the unique positive solution of (\ref{eq 08-10-13 2a})
satisfying
\begin{equation}\label{zero-1}
\lim_{t\to 0}\myint{B_\epsilon}{}u(t,x) dx=\infty,\qquad\forall\epsilon >0,\;B_\epsilon=B_\epsilon(0),
\end{equation}
and
\begin{equation}\label{zero-2}
\lim_{t\to 0}\myint{K}{}u(t,x) dx=0\qquad,\forall K\subset\R^N\setminus\{0\},\, K\text{ compact}.
\end{equation}
If one replaces $u^p$ by $t^\beta u^p$ with $p\in(1,1+\frac{2(1+\beta)}N)$, these  results were extended by Marcus and V\'{e}ron ($\beta\geq 0$) in \cite {MV3} and then Al Sayed and V\'{e}ron ($\beta>-1$) in \cite{SV0}. The initial data problem with measure and general absorption term
\begin{equation}\label{he 1.2}
 \arraycolsep=1pt
\begin{array}{lll}
\partial_t u -\Delta  u+h(t,x,u)=0\quad & \rm{in}\quad (0,T)\times \Omega,\\[2mm]
 \phantom{\partial_t + -\Delta  +h(t,u)= }
u=0\quad&  {\rm in}\quad(0,T)\times\partial\Omega,\\[2mm]
 \phantom{\partial_t u -\Delta u + \ \ =  }
u(0,\cdot)=\nu\quad & \rm{in}\quad \Omega,
\end{array}
\end{equation}
in a bounded domain $\Omega$ of $\R^N$, has been studied by Marcus and V\'{e}ron in \cite {MV3} in the framework of the initial trace theory. They  proved that the following general integrability condition on $h$
\begin{equation}\label{zero-3}
 \arraycolsep=1pt
\begin{array}{lll}
0\leq \abs{h(t,x,r)}\leq \tilde h(t)f(|r|)\qquad,\forall (x,t,r)\in \Omega\times \R_+\times\R\\[2mm]
\myint{0}{T}\tilde h(t)f(\sigma t^{\frac{N}{2}})t^{-\frac{N}{2}}dt<\infty\qquad,\forall\sigma >0\\[3mm]
\text{either }\tilde h(t)=t^\alpha\text{ with }\alpha\geq 0\text{ or }f\text{ is convex, }
\end{array}
\end{equation}
in order that the problem has a unique solution for any bounded measure.
In the particular case with $h(t,x,r)=t^\beta |u|^{p-1}u$, it is fulfilled if  $1<p<1+\frac{2(1+\beta)}N$ and $\beta>-1$, and the very singular solution exists in this range of values.\smallskip

Motivated by a growing number of applications in physics and by important links on the theory of L\'{e}vy process,
 semilinear fractional equations has been attracted much interest in last few years,
   (see e.g. \cite{CCV,CF,CD,CFQ,CV1,CT,DI,FW}).
Recently, in \cite{CV2} we
obtained the existence and uniqueness of a weak solution to semilinear fractional elliptic equation
\begin{equation}\label{1.22a2}
 \arraycolsep=1pt
\begin{array}{lll}
 (-\Delta)^\alpha  u+f(u)=\nu\quad & \rm{in}\quad\Omega,\\[2mm]
 \phantom{   (-\Delta)^\alpha  +g(u)}
u=0\quad & \rm{in}\quad \Omega^c,
\end{array}
\end{equation}
when $\nu$ is a Radon measure and $f$ satisfies a subcritical integrability condition. In \cite{CV1} we studied the the different types of isolated singularities when $f(u)=u^p$ where $1<p<\frac{N}{N-2\alpha}$. In particular, assuming that $0\in\Omega$, we proved that the sequence of solutions $\{u_k\}$ ($k\in\N$) of 
(\ref {1.22a2}), with $\nu=k\delta_0$ converges to infinity when $k\to\infty$, if $p\in (0,1+\frac{2\alpha}{N})$ and it converges to a  solution with a strong singularity at $0$ if $p\in (1+\frac{2\alpha}{N},\frac{N}{N-2\alpha})$.
\smallskip

One purpose of this paper is to study the existence and uniqueness of weak solutions to semilinear fractional heat equation
(\ref{he 1.1}) in a measure framework.
We  first make precise the notion of  weak solution of (\ref{he 1.1})
that we will use in this note.
\begin{definition}\label{weak heat solution}
We say that $u$ is a weak solution of (\ref{he 1.1}), if for any $T>0$, $u\in
L^1(Q_T)$, $h(t,u)\in L^1(Q_T)$ and
\begin{equation}\label{weak heat sense}
 \arraycolsep=1pt
\begin{array}{lll}
\myint{Q_T}{} \left(u(t,x)[-\partial_t\xi(t,x)+(-\Delta)^\alpha\xi(t,x)]+h(t,u)\xi(t,x)\right) dx dt\\[2mm]
\phantom{}
=\myint{\R^N}{} \xi(0,x)d\nu-\myint{\R^N}{}\xi(T,x)u(T,x)dx\qquad \forall \xi\in \mathbb{Y}_{\alpha,T},
\end{array}
\end{equation}
where $Q_T=(0,T)\times\R^N$ and $\mathbb{Y}_{\alpha,T}$ is a space of functions $\xi:[0,T]\times \R^N\to\R$ satisfying
\begin{itemize}
\item[]
\begin{enumerate}\item[$(i)$]
$\norm{\xi}_{L^1(Q_T)}+\norm{\xi}_{L^\infty(Q_T)}+\norm{\partial_t\xi}_{L^\infty(Q_T)}+\norm{(-\Delta)^\alpha \xi}_{L^\infty(Q_T)}<+\infty;$
\end{enumerate}
 \begin{enumerate}\item[$(ii)$]
 for $t\in(0,T)$, there exist $M>0$ and $\epsilon_0>0$ such that for all $\epsilon\in(0,\epsilon_0]$,

 $\norm{(-\Delta)_\epsilon^\alpha \xi(t,\cdot)}_{L^\infty(\R^N)}\leq M.$
\end{enumerate}
\end{itemize}
\end{definition}

Before stating our  main theorems, we introduce the subcritical integrability condition for the nonlinearity $h$, that is,
\begin{itemize}
\item[$(H)\ $]
\begin{enumerate}\item[$(i)$]
The function $h:(0,\infty)\times\R\to\R $ is continuous and for any $t\in(0,\infty)$,  $h(t,0)=0$ and $h(t,r_1)\ge h(t,r_2)$ if $r_1\ge r_2$.
\end{enumerate}
 \begin{enumerate}\item[$(ii)$]
There exist $\beta>-1$ and  a continuous, nondecreasing function $g:\R_+\to\R_+$
 such that
$$|h(t,r)|\le t^\beta g(|r|)\qquad \forall (t,r)\in (0,\infty)\times \R$$
and
\begin{equation}\label{1.4}
\arraycolsep=1pt
 \int_1^{+\infty} g(s)s^{-1-p^*_{\beta}}ds<+\infty,
\end{equation}
where
\begin{equation}\label{1.5}
p^*_{\beta}=1+\frac {2\alpha(1+\beta)}N.
\end{equation}
\end{enumerate}

\end{itemize}

We denote by $H_\alpha:(0,\infty)\times\R^N\times\R^N\to \R_+$ the heat kernel for $(-\Delta)^{\alpha}$ in $(0,\infty)\times\R^N$,
by  $\mathbb{H}_\alpha[\nu]$ the associated heat potential of $\nu\in \mathfrak{M}^b(\R^N)$, defined by
$$\mathbb{H}_\alpha[\nu](t,x)=\int_{\R^N} H_\alpha(t,x,y)d\nu(y)$$
and by $\mathcal{H}_\alpha[\mu]$ the Duhamel operator defined for $(t,x)\in Q_T$ and any $\mu\in L^1(Q_T)$ by
$$\mathcal{H}_\alpha[\mu](t,x)=\myint{0}{t}\mathbb{H}_\alpha[\mu(s,.)](t-s,x)ds=\int_0^t\!\!\int_{\R^N} H_\alpha(t-s,x,y)\mu(s,y)dyds.
$$
Now we state our first theorem as follows.
\begin{teo}\label{teo 1}
Assume that  $\nu\in\mathfrak{M}^b(\R^N)$ and  the function $h$
 satisfies  $(H)$.
Then  problem
(\ref{he 1.1}) admits a unique weak solution $u_{\nu}$ such that
\begin{equation}\label{12-09-2}
\mathbb{H}_\alpha[\nu]-\mathcal{H}_\alpha[h(.,\mathbb{H}_\alpha[\nu_+])]\le u_\nu
\le \mathbb{H}_\alpha[\nu]-\mathcal{H}_\alpha[h(.,-\mathbb{H}_\alpha[\nu_-])]\quad {\rm in}\  Q_\infty,
\end{equation}
where  $\nu_+$ and $\nu_-$ are respectively the positive and negative part in the Jordan
decomposition of $\nu$.
Furthermore,
 \begin{itemize}
\item[$(i)\ $] if $\nu$ is nonnegative, so is $u_\nu$;
\item[$\ \ (ii) $] the
mapping: $\nu\mapsto u_\nu$ is increasing and stable in the sense that if $\{\nu_n\}$ is a sequence of positive bounded Radon measures converging to $\nu$ in the weak sense of measures, then $\{u_{\nu_n}\}$ converges to $u_\nu$ locally uniformly in $Q_\infty$.
\end{itemize}

\end{teo}

According to Theorem \ref{teo 1}, there exists a unique positive weak solution $u_k$ to
\begin{equation}\label{2.1}
 \arraycolsep=1pt
\begin{array}{lll}
\partial_t u + (-\Delta)^\alpha  u+ t^\beta u^p=0\quad & {\rm in}\quad Q_\infty,\\[2mm]
\phantom{ \partial_t u + -\Delta^\alpha  t^\beta u^p }
u(0,\cdot)=k\delta_{0}\quad & {\rm in}\quad \R^N,
\end{array}
\end{equation}
where $\beta>-1$, $k>0$ and $p\in(0,p_\beta^*)$.  We observe that  $u_k\to\infty$ in $(0,\infty)\times\R^N$ as $k\to\infty$ for $p\in(0,1]$, see Proposition \ref{prop4.2} for details.
Our next interest in this paper is to study the limit of $u_k$  as $k\to\infty$ for $p\in(1,p_\beta^*)$, which exists since $\{u_k\}_k$ is an increasing sequence of functions,
 bounded by $\left(\frac{1+\beta}{p-1}\right)^{\frac1{p-1}} t^{-\frac{1+\beta}{p-1}}$, and we set
\begin{equation}\label{13-09-080}
  u_\infty=\lim_{k\to\infty} u_k\quad {\rm in}\quad Q_\infty.
\end{equation}
Actually, $u_\infty$ and  $\{u_k\}_{k}$  are classical solutions to equation
\begin{equation}\label{he 2.1}
 \partial_t u + (-\Delta)^\alpha  u+t^\beta u^p=0\quad  {\rm in}\quad\ Q_\infty,
\end{equation}
see Proposition \ref{re 12-09} for details.

\begin{definition}\label{del 1}

$(i)$ A solution $u$ of (\ref{he 2.1}) is called a self-similar solution if
$$
u(t,x)=t^{-\frac{1+\beta}{p-1}}u(1,t^{-\frac1{2\alpha}}x)\qquad (t,x)\in Q_\infty.
$$
$(ii)$ A solution $u$ of (\ref{he 2.1}) is called  a very singular solution if it vanishes on $\R^N\setminus\{0\}$ at $t=0$ and
$$\lim_{t\to0^+}\frac{u(t,0)}{\Gamma_\alpha(t,0)}=+\infty,$$
where $\Gamma_\alpha :=\mathbb H_\alpha[\delta_0]$ is the fundamental solution of
\begin{equation}\label{07-10-0}
\arraycolsep=1pt
\begin{array}{lll}
\partial_t u+(-\Delta)^\alpha u=0\quad &{\rm in}\quad Q_\infty,\\[2mm]
 \phantom{  \partial_t u \Delta^\alpha+  }
 u(0,\cdot)=\delta_0\quad & {\rm in}\quad\R^N.
\end{array}
\end{equation}
\end{definition}

We remark that for $p\in(1,p_\beta^*)$, a self-similar solution $u$ of (\ref{he 2.1})
is also a very singular solution, since
\begin{equation}\label{17-10-0}
 \lim_{t\to0^+}\Gamma_\alpha(t,0)t^{\frac N{2\alpha}}=c_2,
\end{equation}
for some $c_2>0$. For any  self-similar solution $u$ of (\ref{he 2.1}),
 $v(\eta):=u(1,t^{-\frac1{2\alpha}}x)$ with $\eta=t^{-\frac1{2\alpha}}x$  is a solution of {\it the self-similar equation}
\begin{equation}\label{13-09-07}
(-\Delta)^\alpha v-\frac1{2\alpha}\nabla v\cdot\eta-\frac{1+\beta}{p-1}v+v^p=0\quad{\rm in } \quad \R^N.
\end{equation}
Since $\left(\frac{1+\beta}{p-1}\right)^{\frac1{p-1}}$ is a constant nonzero solution of (\ref{13-09-07}),  the function
\begin{equation}\label{20-10-0}
  U_p(t):=\left(\frac{1+\beta}{p-1}\right)^{\frac1{p-1}}t^{-\frac{1+\beta}{p-1}}\qquad t>0
\end{equation}
 is a flat self-similar solution of (\ref{he 2.1}). It is actually the maximal solution of the ODE $y'+t^\beta y^p=0$ defined on $\R_+$.
Our next goal in this paper is
to study  non-flat self-similar solutions  of (\ref{he 2.1}).

\begin{teo}\label{teo 3}
Assume that $\beta>-1$, $u_\infty$ is defined by (\ref{13-09-080}) and
$$
p_\beta^{**}<p<p_\beta^*,
$$
where $p_\beta^{**}=1+\frac{2\alpha({1+\beta})}{N+2\alpha}$.
 Then  $u_\infty$ is a very singular self-similar solution
of (\ref{he 2.1}) in $Q_\infty$.
Moreover, there exists $c_3>1$ such that
\begin{equation}\label{13-09-08}
 \frac{c_3^{-1}}{1+|x|^{N+2\alpha}}\le u_\infty(1,x)
 \le \frac{c_3\ln(2+|x|)}{1+|x|^{N+2\alpha}}\qquad x\in \R^N.
\end{equation}

\end{teo}

When $p_\beta^{**}<p<p_\beta^*$ with $\beta>-1$, we observe that
$u_\infty$ and $U_p$ are self-similar solutions of (\ref{he 2.1}) and $u_\infty$ is non-flat. Now we are ready to consider
 the uniqueness of non-flat self-similar solution of (\ref{he 2.1}) with decay at infinity,  precisely,  we study the uniqueness of self-similar solution to
\begin{equation}\label{5.1}
\arraycolsep=1pt
\begin{array}{lll}
\partial_t u + (-\Delta)^\alpha  u+ t^\beta u^p=0\quad  &{\rm in}\quad Q_\infty,
\\[2mm]\phantom{ \partial_t u +  }
\lim_{|x|\to\infty} u(1,x)=0.
\end{array}
\end{equation}

We remark that if $u$ is self-similar, then the assumption $\lim_{|x|\to\infty} u(1,x)=0$ is equivalent to $\lim_{|x|\to\infty} u(t,x)=0$ for any $t>0$.
Finally, we state the properties of $u_\infty$ when $1<p\leq p_\beta^{**}$  as follows.

\begin{teo}\label{teo 4}
(i) Assume $1<p<p_\beta^{**}$ and $u_\infty$ is defined by (\ref{13-09-080}).
Then  $u_\infty=U_p$,
where $U_p$ is given by (\ref{20-10-0}).\smallskip

\noindent (ii) Assume $p=p_\beta^{**}$ and $u_\infty$ is defined by (\ref{13-09-080}).
Then  $u_\infty$ is a self-similar solution
of (\ref{he 2.1}) such that
\begin{equation}\label{08-10-2}
 u_\infty(t,x)\ge  \frac{c_4t^{-\frac{N+2\alpha}{2\alpha} }}{1+|t^{-\frac1{2\alpha}}x|^{N+2\alpha}}\qquad (t,x)\in(0,1)\times\R^N,
\end{equation}
for some $c_4>0$.
\end{teo}

We note that Theorem \ref{teo 4}  indicates that there exists no self-similar  solution of  (\ref{he 2.1}) with an initial data $u(0,\cdot)$ vanishing in $\R^N\setminus\{0\}$ if $p\in(1,p_\beta^{**})$,
since $u_\infty$ is the least self-similar solution.
In Theorem \ref{teo 4} part $(ii)$,
we do not know if the self-similar solution is flat or not. From the above theorems, we have the following result.

\begin{teo}\label{teo 5}
(i) Assume $p_\beta^{**}<p<p_\beta^*$. Then problem (\ref{13-09-07}) admits a minimal positive solution $v_\infty$ satisfying
\begin{equation}\label{decay}
\lim_{|\eta|\to\infty}|\eta|^{\frac{2\alpha(1+\beta)}{p-1}}v_\infty(\eta)=0.
\end{equation}
Furthermore,
\begin{equation}\label{decay2}
\frac{c_3^{-1}}{1+|\eta|^{N+2\alpha}}\leq v_\infty(\eta)\leq \frac{c_3\ln(2+|\eta|)}{1+|\eta|^{N+2\alpha}}\qquad \forall\eta\in \R^N
\end{equation}
\noindent (ii) Assume $1<p<p_\beta^{**}$. Then problem (\ref{13-09-07}) admits no positive solution satisfying (\ref{decay}).
\end{teo}

The question of uniqueness of the very singular solution in the case $p_\beta^{**}<p<p_\beta^*$ remains an open problem.\medskip

It is worth comparing the above theorems with the results obtained by Nguyen and V\'eron \cite{NV} concerning the limit, when $k\to\infty$ of the solutions $u=u_k$ of 
\begin{equation}\label{NgVe1}
\arraycolsep=1pt
\begin{array}{lll}
\partial_t u -\Delta  u+ u(\ln (u+1)))^\alpha=0\quad  &{\rm in}\quad Q_\infty,
\\[2mm]\phantom{ \partial_t u -\Delta  u+ u---,}
u(0,.)=k\delta_0\quad  &{\rm in}\quad \R^N,
\end{array}
\end{equation}
where $\alpha>0$. Note that $u_k>0$ and the sequence $\{u_k\}$ is increasing. In this problem, they proved that the diffusion is dominating if  $0<\alpha\leq 1$ and the limit of the $u_k$ is infinite.  If 
$1<\alpha\leq 2$ the absorption dominates, but the limit of the $u_k$ is the maximal solution of the associated ODE, 
$y' + y(\ln (y+1)))^\alpha=0$ on $\R_+$. Finally, if  $\alpha> 2$ the limit of the $u_k$ is a solution with a strong isolated singularity at $(0,0)$, which could be called a very singular solution, although it is not self-similar.\medskip

This paper is organized as follows. In Section 2 we introduce  some properties of
Marcinkiewicz spaces
and Kato's type inequality for  non-homogeneous problems. In Section 3 we prove  Theorem \ref{teo 1}. Section 4 is devoted to investigate the properties of   solutions to
(\ref{2.1}). In Section 5 we give the proof of Theorem \ref{teo 3} and Theorem \ref{teo 4}. Finally, we prove Theorem \ref{teo 5}.

\setcounter{equation}{0}

\section{Linear estimates}

\subsection{The Marcinkiewicz spaces}

We recall the definition and basic properties of the Marcinkiewicz
spaces.

\begin{definition}
Let $\Theta\subset \R^{N+1}$ be an open domain and $\mu$ be a positive
Borel measure in $\Theta$. For $\kappa>1$,
$\kappa'=\kappa/(\kappa-1)$ and $u\in L^1_{loc}(\Theta,d\mu)$, we
set
\begin{equation}\label{mod M}
\|u\|_{M^\kappa(\Theta,d\mu)}=\inf\left\{c\in[0,\infty]:\int_E|u|d\mu\le
c\left(\int_Ed\mu\right)^{\frac1{\kappa'}},\ \forall E\subset \Theta,\, \text{E Borel set}\right\}
\end{equation}
and
\begin{equation}\label{spa M}
M^\kappa(\Theta,d\mu)=\{u\in
L_{loc}^1(\Theta,d\mu):\|u\|_{M^\kappa(\Theta,d\mu)}<\infty\}.
\end{equation}
\end{definition}

$M^\kappa(\Theta,d\mu)$ is called the Marcinkiewicz space of
exponent $\kappa$ or weak $L^\kappa$ space and
$\|.\|_{M^\kappa(\Theta,d\mu)}$ is a quasi-norm. The following
property holds.

\begin{proposition}\label{pr 0} \cite{BBC,CV2}
Assume that $1\le q< \kappa<\infty$ and $u\in L^1_{loc}(\Theta,d\mu)$.
Then there exists  $c_5>0$ dependent of $q,\kappa$ such that
$$\int_E |u|^q d\mu\le c_5\|u\|_{M^\kappa(\Theta,d\mu)}\left(\int_E d\mu\right)^{1-q/\kappa},$$
for any Borel set $E$ of $\Theta$.
\end{proposition}

\bremark If $\Omega$  is a smooth domain of  $\R^N$,
we denote by $H^\Omega_\alpha:(0,\infty)\times\Omega\times\Omega\to \R_+$ the heat kernel for $(-\Delta)^\alpha$
and, if $\nu\in \mathfrak{M}^b(\Omega)$, by  $\mathbb{H}^\Omega_\alpha[\nu]$ the corresponding heat potential of $\nu$ defined by
$$\mathbb{H}^\Omega_\alpha[\nu](t,x)=\int_\Omega H^\Omega_\alpha(t,x,y)d\nu(y).$$
When $\Omega=\R^N$, by Fourier transform, it is clear that
$$H_\alpha(t,x,y)=\frac1{(2\pi)^{N/2}}\int_{\R^N}e^{i(x-y)\cdot\zeta- t|\zeta|^{2\alpha}}d\zeta=H_\alpha(t,x-y,0).$$
Furthermore, $\norm {H_\alpha(t,.,0)}_{L^1}$ is independent of $t$. This implies
\begin{equation}\label{markov}
\norm {\mathbb{H}^\Omega_\alpha[\nu](t,.)}_{L^p}\leq \norm {\nu}_{L^p},\qquad\forall 1\leq p\leq\infty\,,\;\forall \nu\in L^p(\R^N).
\end{equation}
Since $\mathbb{H}^\Omega_\alpha[\nu](t+s,.)=\mathbb{H}^\Omega_\alpha[\mathbb{H}^\Omega_\alpha[\nu](s,.)](t,.)$ for all $t,s>0$ (semigroup property) and $\nu\geq 0\Longrightarrow \mathbb{H}^\Omega_\alpha[\nu](t,.)\geq 0$ the semigroup $\{\mathbb{H}^\Omega_\alpha[.](t,.)\}_{t\geq 0}$ is sub-Markovian. Furthermore, since  the operator $(-\Delta)^\alpha$ is symmetric in $L^2(\R^N)$, the above semigroup is analytic in $L^p(\R^N)$ for all $1\leq p<\infty$: if $1<p<\infty$ it follows from a general result of Stein \cite{St}) and for $p=1$ it is a consequence of regularity result from fractional powers of operators theory (see e.g. \cite{Ko}). For $1\leq p<\infty$ the generator  $A_p$ of the semigroup in $L^p(\R^N)$ is the operator $-(-\Delta)^{\alpha}$ with domain
\begin{equation}\label{markov1}
D(A_p):=\{\nu\in L^p(\R^N): (-\Delta)^{\alpha}\nu\in L^p(\R^N)\}.
\end{equation}
and $D(A_p)$ is dense since it contains $C^\infty_0(\R^N)$. If $p=\infty$, the natural space is the space $C_0(\R^N)$ of continuous functions in $\R^N$ tending to $0$ at infinity. The domain of the corresponding operator $A_{c_0}$ is
\begin{equation}\label{markov2}
D(A_{c_0}):=\{\nu\in C_0(\R^N): (-\Delta)^{\alpha}\nu\in C_0(\R^N)\}.
\end{equation}
This operator is densely defined in $C_0(\R^N)$. In order to avoid confusion, $C_c(\R^N)$ (resp. $C_c^\infty(\R^N)$) denotes the space of continuous (resp. $C^\infty$) functions in $\R^N$ with compact support. It is a dense subset of $C_0(\R^N)$.
\end{remark}

The following regularizing effect $L^p(\R^N)\mapsto L^q(\R^N)$ ($1\leq p\leq q\leq \infty$) is valid for any submarkovian semigroup of contractions in all $L^p(\R^N)$-spaces which has a self-adjoint generator in $L^2(\R^N)$ (see e.g. \cite{St-0}).
\begin{proposition}\label{pr 01} Assume  $1\leq p\leq q\leq \infty$, $p\neq\infty$. Then for any $\nu\in L^p(\R^N)$, $\mathbb{H}_\alpha[\nu](t,.)\in L^q(\R^N)\cap D(A_q)$ for all $t>0$ and there holds, for some positive constant $c=c(\alpha,N,p,q)$,
\begin{equation}\label{markov3}
\norm{\mathbb{H}_\alpha[\nu](t,.)}_{L^q(\R^N)}\leq \myfrac{c}{t^{\frac{N}{2\alpha}(\frac{1}{p}-\frac{1}{q})}}
\norm{\nu}_{L^q(\R^N)}.
\end{equation}
\end{proposition}

Note also that the function $(t,x)\mapsto \mathbb{H}_\alpha[\nu](t,x)$ is $C^\infty$ in $Q_\infty$ as a result of the analyticity on the semigroup $\{\mathbb{H}_\alpha[.](t)\}_{t>0}$.

\begin{proposition}\label{pr 1}
For any $\beta>-1$ and $T>0$, there exists  $c_6>0$ dependent of $N,\alpha,\beta$ such that
for $\nu \in\mathfrak{M}^b(\Omega)$,
\begin{equation}\label{2.6}
\|\mathbb{H}^\Omega_\alpha[|\nu|]\|_{M^{p^*_{\beta}}(Q^\Omega_T, t^\beta dxdt)}\le c_6\|\nu\|_{\mathfrak{M}^b(\Omega)},
\end{equation}
where $p^*_{\beta}$ is defined by (\ref{1.5}) and $Q^\Omega_T=(0,T)\times \Omega$.

\end{proposition}

In order to prove this proposition, we  introduce some notations.
 For $\lambda>0$ and $y\in\Omega$, let us denote
$$A^\Omega_\lambda(y)=\{(t,x)\in Q^\Omega_T: H^\Omega_\alpha(t,x,y)>\lambda\}\,
\text{ and }\,
 m^\Omega_\lambda(y)=\int_{A^\Omega_\lambda(y)}  t^{\beta}dxdt.$$
We also set $A^{\R^N}_\lambda=A_\lambda$ and $ m^{\R^N}_\lambda= m_\lambda$.

\begin{lemma}\label{heat 1}
There exists $c_7>0$ such that for any $\lambda>1$,
\begin{equation}\label{heat 2}
A_\lambda(y)\subset(0,\ c_7\lambda^{-\frac{2\alpha}N}]\times B_{c_7\lambda^{-\frac1N}}(y),
\end{equation}
where $B_r(y)$ is the ball with radius $r$ and center $y$ in $\R^N$.
\end{lemma}
{\bf Proof.}
We observe that $ H_\alpha(t,x,y)=t^{-\frac N{2\alpha}}\Gamma_\alpha (1,(x-y)t^{-\frac1{2\alpha}})$,
where $\Gamma_\alpha $ is the fundamental solution of (\ref{07-10-0}).
From  \cite{BG} (see also\cite{CKS} for an analytic proof), there exists $c_8>0$ such that
$$\Gamma_\alpha (1,z)\le\frac{c_8}{1+|z|^{N+2\alpha}}.$$
This implies in particular
\begin{equation}\label{H-alfa}
H_\alpha(t,x,y)\leq \myfrac{{c_8}t^{-\frac{N}{2\alpha}}}{1+\left(t^{-\frac{1}{2\alpha}}|x-y|\right)^{N+2\alpha}}.
\end{equation}
On the one hand, for $(t,x)\in A_\lambda(y)$, we have that
$$t^{-\frac N{2\alpha}}\Gamma_\alpha(1,0)\ge t^{-\frac N{2\alpha}}\Gamma_\alpha(1,(x-y)t^{-\frac1{2\alpha}})>\lambda,$$
which implies
\begin{equation}\label{2.8}
t<\Gamma_\alpha^{\frac{2\alpha}N}(1,0)\lambda^{-\frac{2\alpha}N}.
\end{equation}
On the other hand, letting $r=|x-y|$,
$$\frac{c_8t}{t^{1+\frac{N}{2\alpha}}+r^{N+2\alpha}} \ge t^{-\frac N{2\alpha}}\Gamma_\alpha(1,(x-y)t^{-\frac1{2\alpha}})>\lambda,$$
then
\begin{equation}\label{2.7}
r\le (c_8t\lambda^{-1})^{\frac1{N+2\alpha}},
\end{equation}
which, together with (\ref{2.8}), implies
$$r\le c_9\lambda^{-\frac1N},$$
for some $c_9>0$.
\hfill$\Box$\\

\noindent{\bf Proof of Proposition \ref{pr 1}.} By Lemma \ref{heat 1}, there exists $c_{10}>0$ such that
$$m_\lambda(y)\le c_{10}\lambda^{-1-\frac{2\alpha(1+\beta)}N}.$$
Clearly
\begin{equation}\label{2.02}
H^\Omega_\alpha(t,x,y)\le H_\alpha(t,x,y),
\end{equation}
then for any Borel set $E\subset Q^\Omega_T$ and $y\in\Omega$, we have that
\begin{eqnarray*}
\myint{E}{} H^\Omega_\alpha(t,x,y)t^\beta dxdt\le\lambda\myint{E}{}t^{\beta} dxdt+\myint{A_\lambda(y)}{} H_\alpha(t,x,y) t^{\beta} dxdt
\end{eqnarray*}
and
$$
\arraycolsep=1pt
\begin{array}{lll}
\myint{A_\lambda(y)}{} H_\alpha(t,x,y)t^{\beta} dxdt=-\myint{\lambda}{+\infty} s d m_s(y)=\lambda m_\lambda(y)+\myint{\lambda}{+\infty}  m_s(y)ds\\[4mm]
\phantom{\myint{A_\lambda(y)}{} H_\alpha(t,x,y)t^{\beta} dxdt}
\le c_{10}\lambda^{-\frac{2\alpha(1+\beta)}N}+c_{10}\myint{\lambda}{+\infty}  s^{-1-\frac{2\alpha(1+\beta)}N}ds
\\[3mm]\phantom{\myint{A_\lambda(y)}{} H_\alpha(t,x,y)t^{\beta} dxdt}
\le c_{11} \lambda^{-\frac{2\alpha(1+\beta)}N},
\end{array}
$$
where $c_{11}=c_{10}\left(1+\frac N{2\alpha(1+\beta)}\right)$.
As a consequence, it follows
\begin{eqnarray*}
\int_{E} H^\Omega_\alpha(t,x,y)t^\beta dxdt\le \lambda\int_{E}  t^{\beta} dxdt+c_{11} \lambda^{-\frac{2\alpha(1+\beta)}N}.
\end{eqnarray*}
Taking $\lambda=(\int_{E}  t^{\beta} dxdt)^{-\frac N{N+2\alpha(1+\beta)}}$, we obtain that
\begin{eqnarray}\label{2.12}
\int_{E} H^\Omega_\alpha(t,x,y)t^\beta dxdt\le (c_{11}+1)(\int_{E}  t^{\beta} dxdt)^{\frac {2\alpha(1+\beta)}{N+2\alpha(1+\beta)}}.
\end{eqnarray}
Since,  by Fubini's theorem,
\begin{eqnarray*}
\int_{E}\mathbb{H}^\Omega_\alpha[|\nu|](t,x) t^{\beta}dxdt
&=&\int_E\int_\Omega H^\Omega_\alpha(t,x,y)d|\nu(y)|t^\beta dxdt
\\&=&\int_\Omega\int_E H^\Omega_\alpha(t,x,y)t^\beta dxdtd|\nu(y)|,
\end{eqnarray*}
together with (\ref{2.12}), it yields
$$\int_{E}\mathbb{H}^\Omega_\alpha[|\nu|](t,x) t^{\beta}dxdt\le (c_{11}+1)\|\nu\|_{\mathfrak{M}^b(\Omega)}\left(\int_E t^{\beta} dxdt\right)^{\frac {2\alpha(1+\beta)}{N+2\alpha(1+\beta)}}.$$
Thus,
\begin{eqnarray*}
\|\mathbb{H}^\Omega_\alpha[|\nu|]\|_{M^{1+\frac {2\alpha(1+\beta)}N}
(Q^\Omega_T,t^{\beta} dxdt)}\le
(c_{11}+1)\|\nu\|_{\mathfrak{M}^b(\Omega)},
\end{eqnarray*}
which ends the proof.\hfill$\Box$


\subsection{The non-homogeneous problem}

In this section we consider the linear non-homogeneous problem
\begin{equation}\label{homo}
 \arraycolsep=1pt
\begin{array}{lll}
\partial_t u + (-\Delta)^\alpha  u=\mu\quad & {\rm in}\quad Q_T,\\[2mm]
\phantom{\partial   (\Delta)^\alpha\  \  }
u(0,\cdot)=\nu\quad & {\rm in}\quad \R^N.
\end{array}
\end{equation}

If $\mu\in L^1(Q_T)$ and $\nu\in L^1(\R^N)$ a function $u$ defined in $Q_T$ is an {\it integral solution} of (\ref{homo}) in $Q_T$ if it is expressed by Duhamel's formula, that is
\begin{equation}\label{Duham}
 \arraycolsep=1pt
\begin{array}{lll}
u(t,x)=\mathbb{H}_\alpha[\nu](t,x)+\mathcal{H}_\alpha[\mu](t,x)\qquad\text{a.e. in } Q_T.
\end{array}
\end{equation}
where, we denote by $\cal H_\alpha$ the operator of $L^1(Q_T)$ defined for all $(x,t)\in {Q_T}$ by
\begin{equation}\label{Duh}
{\cal H}_\alpha[\mu](x,t)=\int_0^t\mathbb H_\alpha[\mu(.,s)](x,t-s)ds=\int_0^t\!\!\int_{\R^N}H_\alpha(t-s,x,y)\mu(s,y)dy ds.
\end{equation}
Notice that, by Duhamel's formula, there holds
\begin{equation}\label{Duha1}\begin{array}{lll}
  \norm{u(t,\cdot)}_{L^1(\R^N)} \leq\norm{\mu}_{L^1(Q_T)}+\norm{\nu}_{L^1(\R^N)},\quad\forall t\in (0,T),
\end{array}\end{equation}
and 
\begin{equation}\label{Duha1}
 \norm{u}_{L^1(Q_T)}\le T(\norm{\mu}_{L^1(Q_T)}+\norm{\nu}_{L^1(\R^N)}).
\end{equation}
The advantage of this notion of solution is that Duhamel's formula has a meaning as soon as $\mu$ and $\nu$ are integrable in their respective domains of definition. As for any continuous semigroup of bounded linear operators, a strong solution is an integral solution. \medskip

The following proposition is the Kato's type estimate which is essential tool to prove
the uniqueness of solutions to (\ref{he 1.1}). For $T>0$, we denote $Q_T=(0,T)\times\R^N$.

\begin{proposition}\label{pr 2.2}
Assume  $\mu\in L^1(Q_T)$ and $\nu\in L^1(\R^N)$. Then there exists a unique weak
solution $u\in L^1(Q_T)$ to the problem (\ref{homo})
and there exists $c_{12}>0$ such that
\begin{equation}\label{L2}
\displaystyle\int_{Q_T} |u| dxdt\leq  \displaystyle c_{12}\int_{Q_T} |\mu| dxdt+c_{12}\int_{\R^N}|\nu|dx.
\end{equation}
Moreover, for any $\xi\in \mathbb{Y}_{\alpha,T}$, $\xi\ge0$, we have that
\begin{equation}\label{sign}
 \arraycolsep=1pt
\begin{array}{lll}
\myint{Q_T}{}  |u|(-\partial_t\xi+(-\Delta)^\alpha \xi) dxdt+\myint{\R^N}{}|u(T,x)|\xi(T,x)dx
\\[4mm]
\phantom{------------ }
\le \myint{Q_T}{}  \xi {\rm sign}(u)\mu dxdt+\myint{\R^N}{}\xi(0,x)|\nu| dx
\end{array}
\end{equation}
and
\begin{equation}\label{sign+}
 \arraycolsep=1pt
\begin{array}{lll}
\myint{Q_T}{}  u_+(-\partial_t\xi+(-\Delta)^\alpha \xi) dxdt+\myint{\R^N}{}u_+(T,x)\xi(T,x)dx
\\[4mm]
\phantom{------------ }
\le \myint{Q_T}{}  \xi {\rm sign_+}(u)\mu dxdt+\myint{\R^N}{}\xi(0,x)\nu_+ dx.
\end{array}
\end{equation}
\end{proposition}

In order to prove Proposition \ref{pr 2.2}, we introduce the following notations.
We say that $u:Q_T\to\R$ is in $ C_{t,x}^{\sigma,\sigma'}(Q_T)$ for $\sigma,\sigma'\in(0,1)$ if
 $$\norm{u}_{ C_{t,x}^{\sigma,\sigma'}(Q_T)}:=\norm{u}_{L^\infty(Q_T)}
+\sup_{Q_T}\frac{|u(t,x)-u(s,y)|}{|t-s|^\sigma+|x-y|^{\sigma'}}<+\infty$$
and $u\in  C_{t,x}^{1+\sigma,2\alpha+\sigma'}(Q_T)$ if
$$\norm{u}_{ C_{t,x}^{1+\sigma,2\alpha+\sigma'}(Q_T)}:=\norm{u}_{L^\infty(Q_T)}+\norm{\partial_tu}_{C_{t,x}^{\sigma,\sigma'}(Q_T)}+\norm{(-\Delta)^\alpha u}_{C_{t,x}^{\sigma,\sigma'}(Q_T)}<+\infty.$$

\begin{lemma}\label{lm 24-10-0}
Let $\mu\in C^1( Q_T)\cap L^\infty(Q_T)$, $\nu\in  L^\infty(\R^N)$ and $u$ be an integral solution of  problem (\ref{homo}),
 then there exists $\sigma\in(0,1)$ such that $u\in C_{t,x}^{1+\sigma,2\alpha+\sigma}$ in $(\epsilon,T)\times\R^N$ for any $\epsilon\in(0,T)$.
 In particular, if  $\|D^2\nu\|_{L^\infty(\R^N)}+\|(-\Delta)^\alpha\nu\|_{C^{1-\alpha}_x(\R^N)}<\infty$, then $u\in C_{t,x}^{1+\sigma,2\alpha+\sigma}(Q_T)$.
\end{lemma}
{\bf Proof.} {\it Step 1.} When $\|D^2\nu\|_{L^\infty(\R^N)}+\|(-\Delta)^\alpha\nu\|_{C^{1-\alpha}_x(\R^N)}<\infty$, it follows directly by  \cite[$(A.1)$]{CF}  that $u\in C^{1+\sigma,2\alpha+\sigma}_{t,x}(Q_T)$.\smallskip

\noindent{\it Step 2.} When $\nu\in L^\infty(\R^N)$, we use \cite[Theorem 6.1]{CD} to obtain that $u\in C^{\frac{\sigma}{2\alpha},\sigma}_{t,x}(Q_T)$ for some $\sigma>0$. For any $\epsilon\in(0,T)$, let $\eta:[0,T]\to[0,1]$ be a $C^2$ function such that $\eta=0$ in $[0,\frac{\epsilon}4]$ and $\eta=1$ in $[\epsilon,T]$ and $v=\eta u$ in $Q_T$.
Since $\eta$ does not depend on $x$, we obtain that $v$ satifies,
$$\partial_tv+(-\Delta)^\alpha v=\eta\mu+\eta'(t)u,\qquad\forall (t,x)\in Q_T,$$
where $\eta\mu+\eta'(t)u\in C^{\frac{\sigma}{2\alpha},\sigma}_{t,x}(Q_T)$ and $v(0,\cdot)=0$ in $\R^N$,
Then we apply the argument in Step 1 to obtain that $v\in C^{1+\sigma,2\alpha+\sigma}_{t,x}(Q_T)$. Therefore,  $u$ is $C^{1+\sigma,2\alpha+\sigma}_{t,x}$  in $(\epsilon,T)\times\R^N$.
The proof is complete. \hfill$\Box$


\begin{lemma}\label{lm 1}
$(i)$
Let $\mu\in C^1( Q_T)\cap L^\infty(Q_T)$ and $\nu\in C^1(\R^N)\cap L^\infty(\R^N)$, then problem (\ref{homo})
 admits a unique classical solution $u$.\smallskip

\noindent$(ii)$ Let $\mu\in C^1( Q_T)\cap L^\infty(Q_T)\cap L^1(Q_T)$, $\nu\in C^2(\R^N)\cap L^\infty(\R^N)\cap L^1(\R^N)$ and $u$ be the classical solution of (\ref{homo}), then $u$ is $ C_{t,x}^{1+\sigma,2\alpha+\sigma}$ in $(\epsilon,T)\times\R^N$ for any $\epsilon\in(0,T)$ and for any
 $\xi\in \mathbb{Y}_{\alpha,T}$,
 \begin{equation}\label{27-10-0}
 \arraycolsep=1pt
\begin{array}{lll}
\myint{Q_T}{}u(t,x)[-\partial_t\xi(t,x)+(-\Delta)^\alpha \xi(t,x)]dxdt \\[2mm]
\phantom{--- }
= \myint{Q_T}{}\mu(t,x)\xi(t,x) dxdt+\myint{\R^N}{}\xi(0,x)\nu dx-\myint{\R^N}{}\xi(T,x)u(T,x)dx.
\end{array}
\end{equation}
Thus $u$ is a weak solution and it belongs to $\mathbb{Y}_{\alpha,T}$.\smallskip

\noindent$(iii)$ Let $\tilde\mu\in C^1( Q_T)\cap L^\infty(Q_T)$ and $\nu\in C^1(\R^N)\cap L^\infty(\R^N)$, then problem
\begin{equation}\label{2.01-}
 \arraycolsep=1pt
\begin{array}{lll}
-\partial_t w + (-\Delta)^\alpha  w=\tilde\mu\quad & {\rm in}\quad Q_T,\\[2mm]
 \phantom{- (-\Delta)\  u}
w(T,\cdot)=\nu\quad & {\rm in}\quad \R^N
\end{array}
\end{equation}
admits a unique classical solution $w\in C_{t,x}^{1+\sigma,2\alpha+\sigma}(Q_T)$  for some $\sigma \in(0,1)$.
Moreover,
 if $\mu\in C^1( Q_T)\cap L^\infty(Q_T)\cap L^1(Q_T)$ and $\nu\in C^2(\R^N)\cap L^\infty(\R^N)\cap L^1(\R^N)$, then  $\xi$ is a weak solution and it belongs to $\mathbb{Y}_{\alpha,T}$.
\end{lemma}
{\bf Proof.} $(i)$ By \cite[Theorem 3.3, Theorem 6.1]{CD}, if $\mu$ and $\nu$ are continuous and bounded, there exists a unique viscosity solution $u\in C(\overline Q_T)$. The higher regularity is provided by \cite[Theorem 6.1]{CD} which asserts that there exist $\sigma>0$ and a positive constant $c$ depending on $N$, $\tau\in (0,T)$ and $\alpha$ such that for all $(t,x)$ and $(s,y)$ belonging to 
$Q^{B_1}_{T-\tau}$, there holds
\begin{equation}\label{2.01--}\displaystyle\myfrac{\abs{u(t,x)-u(s,y)}}{(|x-y|+|t-s|^{\frac{1}{2\alpha}})^\sigma}\leq c\left(\norm {u}_{L^\infty(Q^{B_2}_T)}+\displaystyle\sup_{0\leq t\leq T}\norm {u(t,.)}_{L^1(\R^N}+\norm\mu_{L^\infty(Q_T)}
\right)
\end{equation}
where $Q^{\Omega}_T=(0,T)\times \Omega$. Thus $u\in C^{\frac{\sigma}{2\alpha},\sigma}_{t,x}(Q_T)$.
By Lemma \ref{lm 24-10-0} the integral solution $u$ belongs to $C^{1+\sigma',2\alpha+\sigma'}_{t,x}$
 in $(\epsilon,T)\times\R^N$ for any $\epsilon\in(0,T)$ and some $\sigma'\in(0,\min\{\frac{\sigma}{2\alpha},\sigma\})$. Then  $u$ is a classical solution of (\ref{homo}) and thus a viscosity solution.\smallskip

\noindent$(ii) $ By the definition of $(-\Delta)^\alpha u$, $u(t,\cdot)\in L^1(\R^N)$ for all $t\in(0,T)$. As in \cite[Appendix A.2]{CF} we have Duhamel formula, thus  $u\in L^1(Q_T)$ and it is an integral solution.
\smallskip

\noindent{\it We claim that   $\norm{(-\Delta)_\epsilon^\alpha u(t,\cdot)}_{L^\infty(\R^N)}$ is uniformly bounded with respect to $\epsilon\in(0,\epsilon_0)$.}
Since $u(t,\cdot)\in C^{2\alpha+\sigma}_x(\R^N)$ for some $\sigma\in(0,\min\{2-2\alpha,1\})$,
then for $x\in\R^N$ and $y\in B_1(0)$, $|u(x+y)+u(x-y)-2u(x)|\le \norm{u(t,\cdot)}_{C^{2\alpha+\sigma}_x(\R^N)}|y|^{2\alpha+\sigma}$. Thus,
\begin{eqnarray*}
\norm{|(-\Delta)_\epsilon^\alpha u(t,\cdot)|}_{L^\infty(\R^N)}
&\le &\sup_{x\in\R^N}\left[\int_{\R^N\setminus B_1(0)} \frac{|u(x+y)-u(x)|}{|y|^{N+2\alpha}}dy\right. \\&&\left.+\frac12 \int_{ B_1(0)\setminus B_\epsilon(0)} \frac{|u(x+y)+u(x-y)-2u(x)|}{|y|^{N+2\alpha}}dy\right]  \\
  &\le& 2\norm{u}_{L^1(\R^N)}+\int_{B_1(0)}|y|^{\sigma-N}dy\norm{u(t,\cdot)}_{C^{2\alpha+\sigma}_x(\R^N)}.
\end{eqnarray*}
Next we claim that
\begin{equation}\label{29-10-0}
\int_{Q_T}\xi(-\Delta)_{\epsilon}^\alpha
udxdt=\int_{Q_T}u(-\Delta)_{\epsilon}^\alpha \xi
dxdt\qquad  \ \ \forall\xi\in \mathbb{Y}_{\alpha,T}.
\end{equation}
Indeed, using the fact that for any $t>0$ there holds
\begin{eqnarray*}
&&\myint{\R^N}{}\int_{\R^N}\frac{[u(t,z)-u(t,x)]\xi(t,x)}{|z-x|^{N+2\alpha}}\chi_\epsilon(|x-z|)dzdx\\
&&\qquad\qquad\qquad\qquad=\myint{\R^N}{}\int_{\R^N}\frac{[u(t,x)-u(t,z)]\xi(t,z)}{|z-x|^{N+2\alpha}}\chi_\epsilon(|x-z|)dzdx,
\end{eqnarray*}
then we have
$$
 \arraycolsep=1pt
\begin{array}{lll}
\myint{\R^N}{}\xi(t,x)(-\Delta)_{\epsilon}^\alpha
u(t,x)dx\\[4mm]
\phantom{-}
=-\myfrac{1}2\myint{\R^N}{}\myint{\R^N}{}\left[\myfrac{(u(t,z)-u(t,x))\xi(t,x)}{|z-x|^{N+2\alpha}}+\myfrac{(u(t,x)-u(t,z))\xi(t,z)}{|z-x|^{N+2\alpha}}\right]\chi_\epsilon(|x-z|)dzdx
\\[4mm]
\phantom{-}
=\myfrac{1}2\myint{\R^N}{}\myint{\R^N}{}\myfrac{[u(t,z)-u(t,x)][\xi(t,z)-\xi(t,x)]}{|z-x|^{N+2\alpha}}\chi_\epsilon(|x-z|)dzdx.
\end{array}
$$
Similarly,
$$
\arraycolsep=1pt
\begin{array}{lll}
\myint{\R^N}{}u(t,x)(-\Delta)_{\epsilon}^\alpha \xi(t,x)
dx=\myfrac12\myint{\R^N}{}\myint{\R^N}{}\myfrac{[u(t,z)-u(t,x)][\xi(t,z)-\xi(t,x)]}{|z-x|^{N+2\alpha}}\chi_\epsilon(|x-z|)dzdx.
\end{array}
$$
Then (\ref{29-10-0}) holds. Since  $u$ is $ C_{t,x}^{1+\sigma,2\alpha+\sigma}$ in $(\epsilon,T)\times\R^N$ for any $\epsilon\in(0,T)$ and $\xi$ belongs to $\mathbb{Y}_{\alpha,T}$,
$(-\Delta)_\epsilon^\alpha \xi(t,\cdot)\to (-\Delta)^\alpha \xi(t,\cdot)$ and
$(-\Delta)_\epsilon^\alpha u(t,\cdot)\to (-\Delta)^\alpha u(t,\cdot)$ as  $\epsilon\to 0$ in $\R^N$ and
$(-\Delta)_\epsilon^\alpha\xi(t,\cdot),\ (-\Delta)_\epsilon^\alpha
u(t,\cdot) \in L^\infty(\R^N)$ and $\xi(t,\cdot), u(t,\cdot)\in L^1(\R^N)$, then it follows by the Dominated Convergence
Theorem that
$$\lim_{\epsilon\to 0^+}\int_{\R^N}\xi(t,x)(-\Delta)_\epsilon^\alpha
u(t,x)dx=\int_{\R^N}\xi(t,x)(-\Delta)^\alpha u(t,x)dx$$ and $$
\lim_{\epsilon\to 0^+}\int_{\R^N}(-\Delta)_\epsilon^\alpha\xi(t,x)
u(t,x)dx=\int_{\R^N}(-\Delta)^\alpha \xi(t,x)u(t,x)dx.$$ Combining this with (\ref{29-10-0}), and letting
$\epsilon\to0^+$, we have that
$$\int_{\R^N}\xi(t,x)(-\Delta)^\alpha u(t,x)dx=\int_{\R^N}(-\Delta)^\alpha \xi(t,x)u(t,x)dx,$$
 integrating  over $[0,T]$ and by (\ref{homo}), we conclude that (\ref{27-10-0}) holds.

 \smallskip

\noindent$(iii)$ {\it End of the proof.}  Let $u$ be the weak solution of problem (\ref{homo}) obtained from (ii) with $\tilde \mu(T-t,.)=\mu(t,.)$   and
 $$w(t,x)=u(T-t,x)\qquad (t,x)\in[0,T]\times \R^N.$$
Then $w$ is a solution of (\ref{2.01-}) and for some $\sigma\in(0,1)$, $w$ is $ C_{t,x}^{1+\sigma,2\alpha+\sigma}(Q_T)$. On the contrary, if $w$ is a solution of (\ref{2.01-}), then $u(t,x)=w(T-t,x)$ for $(t,x)\in[0,T]\times \R^N$
is a solution of (\ref{homo}), then the uniqueness holds since the solution of (\ref{homo}) is unique.
Since $u\in C_{t,x}^{1+\sigma,2\alpha+\sigma}(Q_T)$, then $(-\Delta)^\alpha u(t,\cdot)\in C^{\sigma}_x$ and then
$(-\Delta)^\alpha_\epsilon u(t,\cdot)$ is bounded, which implies $u\in\mathbb{Y}_{\alpha,T}$.
 \hfill$\Box$
\medskip

\noindent{\bf Proof of Proposition \ref{pr 2.2}.}
{\it Uniqueness.}
Let $v\in L^1(Q_T)$ be a weak solution of \begin{equation}\label{L0}
\arraycolsep=1pt
\begin{array}{lll}
\partial_t v + (-\Delta)^\alpha  v=0\quad & {\rm in}\quad Q_T,\\[2mm]
\phantom{ ---\ \  }
v(0,\cdot)=0\quad & {\rm in}\quad \R^N.
\end{array}
\end{equation}  We  claim that $v=0$ a.e. in $Q_T$.

In fact, let $\omega$ be a Borel subset of $Q_T$ and $\eta_{\omega,n}$ be the
solution of
\begin{equation}\label{L00}
\arraycolsep=1pt
\begin{array}{lll}
-\partial_t u + (-\Delta)^\alpha  u=\zeta_n\quad & {\rm in}\quad Q_T,\\[2mm]
\phantom{ ----\ \  }
u(T,\cdot)=0\quad & {\rm in}\quad \R^N,
\end{array}
\end{equation}
where $\zeta_n:\bar Q_T\to[0,1]$ is a
function $C^1_c(Q_T)$ such that
$$\zeta_n\to\chi_\omega\quad {\rm{in}}\ L^\infty(\bar  Q_T)\quad {\rm{as}}\ n\to\infty.$$
Then $\eta_{\omega,n}\in \mathbb{Y}_{\alpha,T}$ by Lemma \ref{lm 1},
 and
$$\displaystyle\int_{ Q_T}v\zeta_n dxdt=0.
$$
Passing to the limit when $n\to\infty$, we derive
$$\displaystyle\int_{\omega}v dxdt=0.$$ This implies $v=0$  a.e.  in $ Q_T$.
\smallskip

\noindent{\it Existence and estimate (\ref{sign})}. For  $\delta>0$,
we define an even convex function  $\phi_\delta$ by
\begin{equation}\label{2.3}
\phi_\delta(t)=\left\{ \arraycolsep=1pt
\begin{array}{lll}
|t|-\frac\delta2\qquad & {\rm if}\quad |t|\ge \delta ,\\[2mm]
\frac{t^2}{2\delta}\qquad & {\rm if} \quad |t|< \delta/2.
\end{array}
\right.
\end{equation}
 Then for any $t,s\in\R$,  $|\phi_\delta'(t)|\le1$, $\phi_\delta(t)\to|t|$ and
$\phi_\delta'(t)\to\rm{sign}(t)$ when $\delta\to0^+$. Moreover,
\begin{equation}\label{2.5}
\phi_\delta(s)-\phi_\delta(t)\ge \phi_\delta'(t) (s-t).
\end{equation}

Let $\{\mu_n\}$, $\{\nu_n\}$ be two sequences of functions in $C^2_0(Q_T)$, $C^2_0(\R^N)$, respectively, such
that
$$\lim_{n\to\infty} \int_{Q_T}|\mu_n-\mu| dxdt=0,\ \ \ \lim_{n\to\infty} \int_{\R^N}|\nu_n-\nu| dx=0.$$
We denote by $u_n$ the corresponding solution  to (\ref{homo}) where $\mu,\nu$  are replaced by $\mu_n, \nu_n$, respectively. By Lemma \ref{lm 24-10-0} and Lemma \ref{lm 1}$(ii)$,
$u_n\in C_{t,x}^{1+\sigma,2\alpha+\sigma}(Q_T)\cap L^1(Q_T)$ and then we use  Lemma 2.3 in \cite{CV2} and Lemma \ref{lm 1} $(ii)$ to obtain that for any $\delta>0$ and $\xi\in\mathbb{Y}_{\alpha,T},\ \xi\ge0$,
$$
 \arraycolsep=1pt
\begin{array}{lll}
\myint{Q_T}{} \phi_\delta(u_n)[-\partial_t\xi+(-\Delta)^\alpha \xi]
dxdt+\myint{\R^N}{} \xi(T,x)\phi_\delta(u_n(T,x))dx
\\[4mm]\phantom{ ------}
=\myint{Q_T}{}
 \xi[\partial_t\phi_\delta(u_n)+(-\Delta)^\alpha\phi_\delta(u_n)] dxdt +\int_{\R^N}\xi(0,x)\phi_\delta(\nu_n)dx
 \\[4mm]\phantom{ ------}
 \le\myint{Q_T}{} \xi\phi_\delta'(u_n)[\partial_t u_n+(-\Delta)^\alpha u_n] dxdt +\myint{\R^N}{}\xi(0,x)\phi_\delta(\nu_n)dx
\\[4mm]\phantom{ ------}= \myint{Q_T}{}  \xi \phi_\delta'(u_n) \mu_ndxdt  +\myint{\R^N}{}\xi(0,x)\phi_\delta(\nu_n)dx.
\end{array}
$$
Letting $\delta\to 0^+$, we obtain
\begin{equation}\label{L1}
 \arraycolsep=1pt
\begin{array}{lll}
\displaystyle\int_{Q_T} |u_n|[-\partial_t\xi+(-\Delta)^\alpha \xi]dxdt+\int_{\R^N} \xi(T,x)|u_n(T,x)|dx
\\[2mm]\phantom{ --------}
\leq\displaystyle \int_{Q_T}  \xi {\rm {sign}}(u_n) \mu_ndxdt+\int_{\R^N}\xi(0,x)|\nu_n|dx.
\end{array}
\end{equation}
Let $\eta_k$ be the solution of
\begin{equation}\label{22-10-0}
\arraycolsep=1pt
\begin{array}{lll}
-\partial_t u +(-\Delta)^\alpha u=\varsigma_k &\quad {\rm in}\quad Q_T,\\[2mm]
\phantom{---- \ \ }
u(T,\cdot)=0&\quad {\rm in}\quad {\R^N},
\end{array}
\end{equation}
where $\varsigma_k:Q_T\to [0,1]$ is a $C^2_0$ function such that $\varsigma_k=1$ in $(0,T)\times B_k(0)$.
From the proof of  Lemma \ref{lm 1}, $\tilde \eta_k(t,x):=\eta_k(T-t,x)$ satisfies with $\tilde\varsigma_k(t,x)=\varsigma_k(T-t,x)$
$$
\arraycolsep=1pt
\begin{array}{lll}
\partial_tu +(-\Delta)^\alpha u=\tilde\varsigma_k&\quad {\rm in}\quad Q_T,\\[2mm]
\phantom{----  }
u(0,\cdot)=0&\quad {\rm in}\quad {\R^N}.
\end{array}
$$
By Lemma \ref{lm 24-10-0},
$\tilde \eta_k\in C_{t,x}^{1+\sigma,2\alpha+\sigma}(Q_T)$ with some $\sigma\in(0,1)$ and
\begin{eqnarray*}
0 \le\tilde \eta_k(t,x)&\le& c_8\int_t^T\int_{\R^N}\frac{(s-t)^{-\frac N{2\alpha}}}{1+|(s-t)^{-\frac1{2\alpha}}(y-x)|^{N+2\alpha}}dyds \\
        &\le &c_8\int_t^T\int_{\R^N}\frac{dz}{1+|z|^{N+2\alpha}}ds\\
        &=& c_{13}(T-t).
\end{eqnarray*}
Taking $\xi=\eta_k$ in (\ref{L1}),  we derive that
$$\int_{Q_T} |u_n|\chi_{(0,T)\times B_k(0)} dxdt\leq  c_{13}T\displaystyle \int_{Q_T} |\mu_n| dxdt+c_{13}T\int_{\R^N}|\nu_n| dx.$$
Then, letting $k\to\infty$, we have
\begin{equation}\label{L22}
\displaystyle\int_{Q_T} |u_n| dxdt\leq  c_{13}T\displaystyle \int_{Q_T} |\mu_n| dxdt+c_{13}T\int_{\R^N}|\nu_n| dx.
\end{equation}
Similarly,
\begin{equation}\label{L3}
\displaystyle\int_{Q_T} |u_n-u_m| dx\leq  c_{13}T\displaystyle
\int_{Q_T} |\mu_n-\mu_m| dxdt+c_{13}T\int_{\R^N} |\nu_n-\nu_m|dx.
\end{equation}
Therefore, $\{u_n\}_n$ is a Cauchy sequence in $L^1(Q_T)$ and its limit $u$
is a weak solution of (\ref{homo}). Letting $n\to\infty$,
(\ref{sign}) and (\ref{L2}) follow by (\ref{L1}) and (\ref{L22}), respectively. The proof of (\ref{sign+}) is similar.
\hfill$\Box$

\bremark Other classes of uniqueness of solutions of the fractional heat equations exist. In \cite{BPSV} it is proved that any positive strong solution $u\in C([0,T)\times\R^N)$ can be represented by the convolution integral defined by
$$u(t,x)=\myint{\R^N}{}P_t(y)u(0,y)dy
$$
where 
$$P_t(x)=\myfrac{1}{t^{\frac{N}{2\alpha}}}e^{it^{-\frac{1}{2\alpha}}x.\xi-|\xi|^{2\alpha}}.$$ 
However the fact that $u\in L^1(Q_T)$ is a part of the definition of strong solution therein. Furthermore, the notion of weak solution used in this paper differs from ours.
\eremark
\noindent


\setcounter{equation}{0}
\section{Proof of Theorem \ref{teo 1}}

If $h(t,.)$ is monotone nondecreasing, for any $\lambda>0$, $I+\lambda h(t,.)$ is an homeorphism of $\R$ and the inverse function
$J_\lambda(t,.)=(I+\lambda h(t,.))^{-1}$ is a contraction. We define the Yosida approximation by
\begin{equation}\label{Yos-1}
h_\lambda(t,.)=\frac{I-J_\lambda(t,.)}{\lambda}.
\end{equation}
The function $h_\lambda(t,.)$ is monotone nondecreasing, vanishes at $0$ as $h$ does it and it is $\frac{1}{\lambda}$-Lipschitz continuous. Furthermore
\begin{equation}\label{Yos-2}
rh_\lambda(t,r)\uparrow rh(t,r)\quad\text{ as }\lambda\to 0,\qquad\forall r\in\R,
\end{equation}
see \cite[Chap 2, Prop. 2.6]{Bre}. If $u$ is a real valued function we will denote by $h\circ u$ and $h_\lambda\circ u$ respectively the functions
$(t,x)\mapsto h(t,u(t,x))$ and $(t,x)\mapsto h_\lambda(t,u(t,x))$.

\begin{lemma}\label{lm 17-09-0}
Assume that $h$ satisfies $(H)$-(i), $\lambda>0$ and $\phi\in L^1(\R^N)$. Then there exists a unique solution $u_\phi$ of
\begin{equation}\label{he 1.1,L}
 \arraycolsep=1pt
\begin{array}{lll}
\partial_t u + (-\Delta)^\alpha  u+h_\lambda\circ u=0\quad & {\rm in}\quad Q_\infty,\\[2mm]
 \phantom{\partial_t u + (-\Delta)^\alpha u_\lambda \ \ }
u(0,\cdot)=\phi\quad & {\rm in}\quad \R^N.
\end{array}
\end{equation}
Moreover,
\begin{equation}\label{14-10-0}
\mathbb{H}_\alpha[\phi]-\mathcal{H}_\alpha[h_\lambda\circ\mathbb{H}_\alpha[\phi_+])]\le u_\phi
\le \mathbb{H}_\alpha[\phi]-\mathcal{H}_\alpha[h_\lambda\circ(-\mathbb{H}_\alpha[\phi_-])]\quad{\rm in}\  Q_{T},
\end{equation}
where $\phi_{\pm}=\max\{0,\pm\phi\}$ and
\begin{equation}\label{14-10-01}
\norm{u_\phi(t,.)-u_\psi(t,.)}_{L^1}\leq \norm{\phi-\psi}_{L^1},\qquad\forall 1\leq p\leq\infty.
\end{equation}

\noindent $(i)$ $u_\phi\ge0$ if $\phi\ge0$ in $\Omega$; \smallskip

\noindent $(ii)$ the mapping $\phi\mapsto u_\phi$ is increasing.

\end{lemma}
{\bf Proof.} Existence is a consequence of  the Cauchy-Lipschitz-Picard theorem (see \cite [Chap 4]{CazHa}): we write (\ref{he 1.1,L}) under the integral form $u=\mathcal T[u]=\mathbb{H}_\alpha[\phi]-\mathcal H_\alpha[h_\lambda\circ u]$, i.e.
\begin{equation}\label{int-1}
\mathcal T[u](t,.)=\mathbb{H}_\alpha[\phi](t,.)-\myint{0}{t}\mathbb{H}_\alpha[h_\lambda\circ u](t-s,.)ds.
\end{equation}
The space
$C([0,\infty);L^1(\R^N))$ endowed with the norm
$$\norm{w}_{C-L^1}=\sup\left\{e^{-k t}\norm{w(t,.)}_{L^1}:t\geq 0\right\}, $$
($k>\lambda^{-1}$), is a Banach space. Since $u\mapsto h_\lambda(t,u)$ is $\frac{1}{\lambda}$-Lipschitz continuous, the mapping $\mathcal T$ is $\frac{1}{\lambda k}$-Lipschitz continuous in $X_{p}$. Thus it admits a unique fixed point $u_\phi$ which is an integral solution of (\ref{he 1.1,L}).
\begin{equation}\label{fix-1}
u_\phi(t,.)=\mathbb{H}_\alpha[\phi](t,.)-\myint{0}{t}\mathbb{H}_\alpha[h_\lambda\circ u_\phi](t-s,.)ds.
\end{equation}
The semigroup $\{\mathbb{H}_\alpha[.](t,.)\}_{t\geq 0}$ is analytic in $L^1(\R^N)$ since it is generated by the fractional power of a closed operator. It follows from the classical regularity theory for analytic semigroups as exposed in \cite[Sec 6]{Gri} that $u_\phi$ is a strong solution of (\ref{he 1.1,L}). Since it is continuous, it is also a weak solution in the sense that
\begin{equation}\label{a-0}
\begin{array}{lll}
\myint{Q_T}{} \left(u_{\phi}[-\partial_t\xi+(-\Delta)^\alpha\xi]+\xi h_\lambda\circ u_{\phi}\right) dx dt\\[2mm]
\phantom{---------}
=\myint{\R^N}{} \xi(0,x)\phi(x)dx-\myint{\R^N}{}\xi(T,x)u_{\phi}(T,x)dx\qquad \forall \xi\in \mathbb{Y}_{\alpha,T}.
\end{array}
\end{equation}\smallskip
If  $\phi_1, \phi_2\in L^1(\R^N)$ and $u_{\phi_j}$ are the corresponding solutions of (\ref{he 1.1,L}), it follows from the positivity of $H_\alpha$ that
$$\begin{array}{lll}
(u_{\phi_2}-u_{\phi_1})_+\leq (\mathcal H_\alpha[h_\lambda \circ u_{\phi_2}-h_\lambda \circ u_{\phi_1}])_+
\leq\myfrac{1}{\lambda}\mathcal H_\alpha[( u_{\phi_2}-u_{\phi_1})_+].
\end{array}
$$
Therefore,
$$\norm{(u_{\phi_2}(t,.)-u_{\phi_1}(t,.))_+}_{L^p}
\leq\myfrac{1}{\lambda}\myint{0}{t}\norm{(u_{\phi_2}(t-s)-u_{\phi_1}(t-s))_+}_{L^p}ds,
$$
and by Gronwall inequality
$$\norm{(u_{\phi_2}(t)-u_{\phi_1}(t))_+}_{L^p}\leq e^{\frac{t}{\lambda}}\norm{(\phi_2-\phi_1)_+}_{L^p}.
$$
This implies (i) and (ii).
As a consequence,
$$-\mathbb H_\alpha[\phi_-]\leq -u_{\phi_-}\leq u_\phi\leq u_{\phi_+}\leq \mathbb H_\alpha[\phi_+]$$
and thus
$$h_\lambda\circ(-\mathbb H_\alpha[\phi_-])\leq h_\lambda\circ(-u_{\phi_-})\leq h_\lambda\circ u_\phi\leq h_\lambda\circ u_{\phi_+}\leq h_\lambda\circ\mathbb H_\alpha[\phi_+].
$$
Jointly with (\ref{fix-1}) it yields (\ref{14-10-0}).
\hfill$\Box$\medskip

\noindent {\bf Notation.} In the sequel, if $\eta\in L^1(Q_\tau)$ and $\tau\geq T$,  we denote by $\xi_{\eta,\tau}$ the solution of
\begin{equation}\label{nota}
\begin{array}{lll}-\partial_t\xi_\eta+(-\Delta)^\alpha\xi_\eta=\eta\qquad\text{in }Q_\tau,\\[2mm]
\phantom{+(-\Delta)^\alpha\xi_\eta}
\xi_\eta(\tau,.)=0.
\end{array}
\end{equation}
If $\eta\geq 0$, then $\xi_{\eta,\tau}\geq 0$; if $\eta\in C^\infty_0(\R^{N+1})$, then $\eta\in\mathbb{Y}_{\alpha,\tau}$; if $\eta_n=\eta(\frac{.}{n})$, where $n\in\mathbb N_*$ and
 $\eta\in C^\infty_0(\R^{N+1})$ is nonnegative, $0\leq \eta\leq 1$, with value $1$ on $B_1$ and $0$ on $B^c_2$, then
$\xi_{\eta_n,\tau}\uparrow \tau-t$ as $n\to\infty$.
\medskip

In the next lemma we prove that we can replace $h_\lambda$ by $h$.

\begin{lemma}\label{lm 2}
Assume  that $h$ satisfies $(H)$-(i) and $\phi\in L^1(\R^N)$. Then there exists a unique solution $u_\phi\in C([0,\infty);L^1(\R^N)$ of
\begin{equation}\label{he 1.1,0}
 \arraycolsep=1pt
\begin{array}{lll}
\partial_t u + (-\Delta)^\alpha  u+h\circ u=0\quad & {\rm in}\quad Q_\infty,\\[2mm]
 \phantom{\partial_t u + (-\Delta)^\alpha u \ \ }
u(0,\cdot)=\phi\quad & {\rm in}\quad \R^N.
\end{array}
\end{equation}
Moreover inequality (\ref{14-10-01}) and statements (i) and (ii) in Lemma \ref{lm 17-09-0} hold.
\end{lemma}
{\bf Proof.} We denote by $u_{\lambda,\phi}$ the solution of (\ref{he 1.1,L}). \smallskip

\noindent{\it Step 1- A priori estimate. } Let $\phi\geq 0$. If we take $\xi=\xi_{\eta_n, \tau}$ in (\ref{a-0}) and let $n\to\infty$, we derive
\begin{equation}\label{a-priori}
\myint{Q_T}{} \left(u_{\lambda,\phi}+(\tau-t)h_\lambda\circ u_{\lambda,\phi}\right) dx dt+(\tau-T)\myint{\R^N}{}u_{\lambda,\phi}(T,x)dx= \tau\myint{\R^N}{}\phi(x)dx.
\end{equation}
For $0<\lambda<\lambda'$ we set
$w=u_{\lambda,\phi}-u_{\lambda',\phi}$. It follows from (\ref{sign+}) and inequality $h_{\lambda'}\circ u_{\lambda,\phi}\leq h_\lambda\circ u_{\lambda,\phi}$, that for any nonnegative $ \xi$ in $\mathbb{Y}_{\alpha,T}$,
$$
\begin{array}{lll}
\myint{Q_T}{} \left(w_{+}[-\partial_t\xi+(-\Delta)^\alpha\xi]+\xi \left({h_\lambda}\circ u_{\lambda,\phi}-h_\lambda\circ u_{\lambda',\phi}\right){\rm sign}_+(w)\right) dx dt\\[4mm]
\phantom{----}
\leq\myint{Q_T}{} w_{+}\left({h_{\lambda'}}\circ u_{\lambda',\phi}-h_\lambda\circ u_{\lambda',\phi}\right)dxdt-\myint{\R^N}{}\xi(T,x)w_{+}(T,x)dx,
\end{array}
$$
Since $h_\lambda(t,.) $ is nondecreasing, we derive
$$\myint{Q_T}{} w_{+}[-\partial_t\xi+(-\Delta)^\alpha\xi]dxdt\leq 0\qquad \forall \xi\in \mathbb{Y}_{\alpha,T},\, \xi\geq 0.
$$
If $\eta\in C^\infty_0(\R^{N+1})$ is nonnegative, then $\xi_\eta\in \mathbb{Y}_{\alpha,T}$, $\xi_\eta\geq 0$ and
$$\myint{Q_T}{} w_{+}\eta dxdt=0.
$$
This implies $u_{\lambda,\phi}\leq u_{\lambda',\phi}$. \smallskip

\noindent{\it Step 2- Truncation.}
We replace $\phi$ by $\phi_n=\inf\{\phi,n\}$ for $n\in\mathbb N_*$ and denote by $u_{\lambda,\phi_n}$ the corresponding solution of (\ref{he 1.1,L}). By Step 1, the sequence $\{u_{\lambda,\phi_n}\}_{\lambda>0}$ is decreasing and it converges to some nonnegative $u_{\phi_n}$ when $\lambda\downarrow 0$. Therefore $h_\lambda\circ u_{\lambda,\phi_n}\to h\circ u_{\phi_n}$ a.e. in $Q_T$. It follows from (\ref{a-priori}) and Fatou's lemma that
\begin{equation}\label{a-priori-n}
\myint{Q_T}{} \left(u_{\phi_n}+(\tau-t)h\circ u_{\phi_n}\right) dx dt+(\tau-T)\myint{\R^N}{} u_{\phi_n}(.,T)dx= \tau\myint{\R^N}{}\phi_n(x)dx.
\end{equation}
Since $0\leq u_{\lambda,\phi_n}\leq n$, then
$0\leq h_\lambda\circ u_{\lambda,\phi_n}\leq h\circ u_{\lambda,\phi_n}\leq h(n)$ by (\ref{14-10-01}). If $E\subset Q_T$ is a Borel set,
$$\myint{E}{}h_\lambda\circ u_{\lambda,\phi_n} dxdt\leq h(n)|E|.
$$
By Vitali convergence theorem $h_\lambda\circ u_{\lambda,\phi_n}\to h\circ u_{\phi_n}$ in $L^1(Q_T)$. Therefore, we can let $\lambdaÊ\to 0$ in identity (\ref{a-0}) and conclude that $u_{\phi_n}$ is a weak solution of (\ref{he 1.1,0}) with initial data $\phi_n$. \smallskip

\noindent{\it Step 3- Existence with $\phi$ bounded. } If $\phi=\phi_+-\phi_-\in L^1(\R^N)$, set $\phi_{+,n}=\inf\{\phi_+,n\}$ and $\phi_{-,n}=\inf\{\phi_-,n\}$. We denote by $u_{\lambda,\phi_{+,n} }$, $u_{\phi_{+,n} }$, $u_{\lambda,-\phi_{-,n} }$ and $u_{-\phi_{-,n} }$ the corresponding solutions of (\ref{he 1.1,L}) and (\ref{he 1.1,0}). Then

\begin{equation}\label{upper} \begin{array} {ccl}
u_{\lambda,-\phi_{-,n} }\leq u_{\lambda,\phi_{+,n}-\phi_{-,n} }\leq u_{\lambda,\phi_{+,n} }\\[2mm]
\text{which implies}\\
h_\lambda\circ u_{\lambda,-\phi_{-,n} }\leq h_\lambda\circ u_{\lambda,\phi_{+,n}-\phi_{-,n} }\leq h_\lambda\circ u_{\lambda,\phi_{+,n} }.
\end{array}\end{equation}
Estimate (\ref{a-priori}) is valid under the form
\begin{equation}\label{a-priori-+}\begin{array} {lll}
\myint{Q_T}{} \left(u_{\lambda,\phi_{+,n}}+(\tau-t)h_\lambda\circ u_{\lambda,\phi_{+,n}}\right) dx dt\\[4mm]
\phantom{-----------}+
(\tau-T)\myint{\R^N}{} u_{\lambda,\phi_{+,n}}(.,T)dx= \tau\myint{\R^N}{}\phi_{+,n}(x)dx.
\end{array}\end{equation}
and
\begin{equation}\label{a-priori-+-}\begin{array} {lll}
\myint{Q_T}{} \left(u_{\lambda,-\phi_{-,n}}+(\tau-t)h_\lambda\circ u_{\lambda,-\phi_{-,n}}\right) dx dt\\[4mm]
\phantom{---------}
+(\tau-T)\myint{\R^N}{} u_{\lambda,-\phi_{-,n}}(.,T)dx= -\tau\myint{\R^N}{}\phi_{-,n}(x)dx.
\end{array}\end{equation}
Since $h_\lambda\circ u_{\lambda,\phi_{+,n}}$ and $h_\lambda\circ u_{\lambda,-\phi_{-,n}}$ are bounded in $L^1(Q_T)$ independently of $\lambda$ and $n$, $h_\lambda\circ u_{\lambda,\phi_{+,n}-\phi_{-,n} }$ inherits the same property. Since
$$u_{\lambda,\phi_{+,n}-\phi_{-,n} }=\mathbb H_\alpha [\phi_{+,n}-\phi_{-,n}]-\mathcal H_\alpha [h_\lambda\circ u_{\lambda,\phi_{+,n}-\phi_{-,n} }],
$$
it follows from \cite[Sec 6]{Gri} that $u_{\lambda,\phi_{+,n}-\phi_{-,n} }$ remains bounded in the interpolation space
$Y_1:=L^1([0,T];D(A_1)(\R^N))\cap W^{s,1}([0,T];L^{1}(\R^N))$, for any $s\in (0,1)$, where $D(A_1)$ is defined in (\ref{markov1}). Although a bounded subset $K$ of $Y_1$ is not a relatively compact subset of $L^1(Q_T)$, for any ball $B\subset\R^N$, the set of restrictions to $B$ of functions belonging to $K$ is relatively compact in $L^1((0,T)\times B)$. Thus, there exists a subsequence
$\{\lambda_k\}$ such that $\{u_{\lambda_k,\phi_{+,n}-\phi_{-,n} }\}$ converges a.e. to some function $U_n$. Furthermore
$\{h_{\lambda_k}\circ u_{\lambda_k,\phi_{+,n}-\phi_{-,n} }\}$ converges a.e. to $h\circ U_n$. Since the sequences
$\{ u_{\lambda_k,-\phi_{-,n} }\}_{\lambda_k}$, $\{ u_{\lambda_k,\phi_{+,n} }\}_{\lambda_k}$,
$\{h_{\lambda_k}\circ u_{\lambda_k,-\phi_{-,n} }\}_{\lambda_k}$ and $\{h_{\lambda_k}\circ u_{\lambda_k,\phi_{+,n} }\}_{\lambda_k}$ are convergent in $L^1(Q_T)$ they are uniformly integrable. Because of (\ref{upper}) the same property is shared by the two sequences $\{u_{\lambda_k,\phi_{+,n}-\phi_{-,n} }\}_{\lambda_k}$ and $\{h_{\lambda_k}\circ u_{\lambda_k,\phi_{+,n}-\phi_{-,n} }\}_{\lambda_k}$. Letting $\lambda_k$ to $0$ in the identity
\begin{equation}\label{a-lambda-k}
\begin{array}{lll}
u_{\lambda_k,\phi_{+,n}-\phi_{-,n} }(t,.)=\mathbb{H}_\alpha[\phi_{+,n}-\phi_{-,n} ](t,.)-
\myint{0}{t}\mathbb{H}_\alpha[h_{\lambda_k}\circ u_{\lambda_k,\phi_{+,n}-\phi_{-,n} }](t-s,.)ds
\end{array}
\end{equation}
yields
\begin{equation}\label{a-n}
\begin{array}{lll}
U_n(t,.)=\mathbb{H}_\alpha[\phi_{+,n}-\phi_{-,n} ](t,.)-
\myint{0}{t}\mathbb{H}_\alpha[h\circ U_n](t-s,.)ds.
\end{array}
\end{equation}
This implies that $U_n$ is an integral solution, thus a weak solution of (\ref{he 1.1,0})  with initial data
$\phi_{+,n}-\phi_{-,n} ={\rm sgn}(\phi)\inf\{n,|\phi|\}$ and then $U_n=u_{\phi_n}$.\smallskip

\noindent{\it Step 4- Existence with $\phi\in L^1(\R^N)$. } By Kato's inequality (\ref{sign}), we obtain that
$$\begin{array}{lll}
\myint{Q_T}{}  \left(|u_{\phi_k}-u_{\phi_m}|(-\partial_t\xi+(-\Delta)^\alpha \xi) +\xi|h\circ u_{\phi_k}-h\circ u_{\phi_m}|\right)dxdt\\[4mm]
\phantom{---- }
+\myint{\R^N}{}|u_{\phi_k}(T,x)-u_{\phi_m}(T,x)|\xi(T,x)dx\le \myint{\R^N}{}\xi(0,x)|\phi_k-\phi_m| dx,
\end{array}
$$
for $m,k\in\mathbb N_*$ and $\xi\in\mathbb Y_{\alpha,T}$, $\xi>0$. Taking $\xi=\xi_{\eta_n,\tau}$ as in (\ref{nota}) and letting $n\to\infty$ yields
\begin{equation}\label{a-07}\begin{array}{lll}
\myint{Q_T}{}  \left(|u_{\phi_k}-u_{\phi_m}| +(\tau-t)|h\circ u_{\phi_k}-h\circ u_{\phi_m}|\right)dxdt\\[4mm]\phantom{-------}
+(\tau-T)\myint{\R^N}{}|u_{\phi_k}(T,.)-u_{\phi_m}(T,.)| dx
\le \tau\myint{\R^N}{}|\phi_k-\phi_m| dx.
\end{array}
\end{equation}
Since $\{\phi_m\}$ is a Cauchy sequence in $L^1(\R^N)$, $\{u_{\phi_m}\}$ and $\{h\circ u_{\phi_m}\}$ are also Cauchy sequences in  $C(0,T;L^1(\R^N))$ and $L^1(Q_T)$ respectively. Set $U=\lim_{m\to\infty}u_{\phi_m}$, then it satisfies
\begin{equation}\label{a-01}
\begin{array}{lll}
\myint{Q_T}{} \left(U[-\partial_t\xi+(-\Delta)^\alpha\xi]+\xi h\circ U\right) dx dt\\[2mm]
\phantom{-------}
=\myint{\R^N}{} \xi(0,x)\phi(x)dx-\myint{\R^N}{}\xi(T,x)U(T,x)dx\qquad \forall \xi\in \mathbb{Y}_{\alpha,T},
\end{array}
\end{equation}
and it is also an integral solution of (\ref{he 1.1,0}). Thus $u_\phi\in C([0,\infty);L^1(\R^N))$. \medskip

Finally, we end the proof of uniqueness which is a consequence of the inequality below
\begin{equation}\label{a-08}\begin{array}{lll}
\myint{Q_T}{}  \left(|U-U'| +(\tau-t)|h\circ U-h\circ U'|\right)dxdt\\[4mm]\phantom{-------}
+(\tau-T)\myint{\R^N}{}|U(T,.)-U'(T,.)| dx
\le \tau\myint{\R^N}{}|\phi-\phi'| dx,
\end{array}
\end{equation}
valid for two solutions $U$ and $U'$ of problem (\ref{he 1.1,0}) with respective initial data $\phi$ and $\phi'$, the proof of which is the same as the one of (\ref{a-07}). Notice also that statement (i) and (ii) as well as inequality (\ref{14-10-01}) follows by the above approximations.\hfill$\Box$
\medskip

\begin{remark} By the same method it can be proved that for any $p\in (1,\infty)$ and $\phi\in L^{p}(\R^N)$ (resp. $\phi\in C_0(\R^N)$) there exists a unique solution $u_\phi\in C([0,\infty);L^{p}(\R^N))$ (resp. $u_\phi\in C([0,\infty); C_0(\R^N))$) solution of
(\ref{he 1.1,0}). Furthermore (\ref{14-10-01}) holds.
\end{remark}

\medskip

\noindent{\bf Proof of Theorem \ref{teo 1}.} \emph{Existence for $\nu\ge0$. }
We  consider a sequence of nonnegative functions $\{\nu_n\}_n\subset C^2_0(\R^N)$ such that
$\nu_{n}\to\nu $ as $n\to\infty$ in the weak sense  of bounded measures, i.e.
\begin{equation}\label{1.3}
\lim_{n\to\infty}\int_{\R^N}\zeta \nu_{n}dx=\int_{ \R^N}\zeta d\nu \qquad \forall  \zeta\in C(\R^N)\cap L^\infty(\R^N).
\end{equation}
 It follows from the
Banach-Steinhaus theorem that $\norm{\nu_{n}}_{\mathfrak M^b
(\R^N)}$ is bounded independently of $n$ and we assume that $\norm{\nu_{n}}_{\mathfrak M^b
(\R^N)}\le 2\norm{\nu}_{\mathfrak M^b
(\R^N)}$. By Lemma \ref{lm 17-09-0}, we denote by $u_{\nu_n}$  the corresponding solution of (\ref{he 1.1,0}) with initial data $\nu_n$. Then
$u_n$ is nonnegative and satisfies that
\begin{equation}\label{10-09-1}
0\le u_{\nu_n}=\mathbb{H}_\alpha[\nu_{n}] -\mathcal H_\alpha [h\circ u_{\nu_n}]\le \mathbb{H}_{\alpha}[\nu_n]\quad{\rm in}\ Q_T.
\end{equation}
Jointly with $(\ref{2.6})$  it implies
\begin{equation}\label{L8}
\norm {u_{\nu_n}}_{M^{p^*_\beta}(Q_T,t^\beta dxdt)}\leq c_5\|\nu\|_{\mathfrak M^b (\R^N)}.
\end{equation}
We have also the following estimates from (\ref{H-alfa}) and (\ref{a-priori-n})
\begin{equation}\label{L8-1}
u_{\nu_n}(t,x)\leq \mathbb{H}_\alpha[\nu_{n}](t,x)\leq 2c_8t^{-\frac{N}{2\alpha}}\|\nu\|_{\mathfrak M^b (\R^N)},\qquad\forall (t,x)\in Q_T
\end{equation}
and
\begin{equation}\label{a-priori-n'}\begin{array}{lll}
\myint{Q_T}{} \left(u_{\nu_n}+(\tau-t)h\circ u_{\nu_n}\right) dx dt+(\tau-T)\myint{\R^N}{} u_{\nu_n}(.,T)dx= \tau\myint{\R^N}{}\nu_n(x)dx\\[4mm]
\phantom{\myint{Q_T}{} \left(u_{\nu_n}+(\tau-t)h\circ u_{\nu_n}\right) dx dt+(\tau-T)\myint{\R^N}{} u_{\nu_n}(.,T)dx}
\leq
2\tau\|\nu\|_{\mathfrak M^b (\R^N)}.
\end{array}\end{equation}
As in the proof of Lemma \ref{lm 2}-Step 3, using the regularizing properties of the semigroup $\mathbb H_\alpha[.](t)$  (see \cite[Sec 6]{Gri}) we infer that there exists a subsequence $\{u_{\nu_{n_k}}\}$ which converges a.e. in $Q_T$ to some function $U$ and $\{h\circ u_{\nu_{n_k}}\}$ converges a.e. to $h\circ U$.\smallskip

 For $\kappa>0$, we denote  $S_\kappa=\{(t,x)\in Q_T:|u_{n_k }(t,x)|>\kappa\}$  and
$\omega(\kappa)=\int_{S_\kappa} t^{\beta}dxdt$. Then for any Borel
set $E\subset Q_T$
$$\begin{array}{lll}
\myint{}{}\!\!\myint{E}{}h\circ u_{\nu_{n_k}} dx dt\leq \myint{}{}\!\!\myint{E\cap\{u_{\nu_{n_k}}\leq \kappa\}}{}h\circ u_{\nu_{n_k}} dx dt+\myint{}{}\!\!\myint{E\cap S_\kappa}{}h\circ u_{\nu_{n_k}} dx dt\\[4mm]
\phantom{\myint{}{}\!\!\myint{E}{}h\circ u_{\nu_{n_k}} dx dt}
\leq g(\kappa) \myint{}{}\!\!\myint{E}{}t^{\beta}dx dt+
\myint{}{}\!\!\myint{\{S_\kappa}{}t^\beta g(u_{\nu_{n_k}}) dx dt\\[4mm]
\phantom{\myint{}{}\!\!\myint{E}{}h\circ u_{\nu_{n_k}} dx dt}
\leq g(\kappa) \myint{}{}\!\!\myint{E}{}t^{\beta}dx dt-\myint{\kappa}{\infty}
g(s)d\omega(s),
\end{array}$$
where
$$\int_{\kappa}^\infty  g(s)d\omega(s)=\lim_{M\to\infty}\int_{\kappa}^M  g(s)d\omega(s).
$$
By (\ref{mod M}) and (\ref{L8}), $\omega(s)\leq c_{14}s^{-p^*_{\beta}}$, thus
$$\displaystyle\begin{array}{lll}
\displaystyle-\int_{\kappa}^M  g(s)d\omega(s) =-\left[g(s)\omega(s)\!\!\!\!\!\!\!\!^{\phantom{\frac{X^X}{X}}}\right]_{s=\kappa}^{s=M}+\int_{\kappa}^M
\omega(s)d g(s)
\\[4mm]\phantom{\displaystyle-\int_{\kappa}^M  g(s)d\omega(s) }\displaystyle
\leq g(\kappa)\omega(\kappa)-g(M)\omega(M)+c_{14}\int_{\kappa}^M s^{-p^*_{\beta}}d g(s)
\\[4mm]\phantom{\displaystyle-\int_{\kappa}^M  g(s)d\omega(s) }\displaystyle
\leq  g(\kappa)\omega(\kappa)- g(M)\omega(M)+
c_{14}\left(M^{-p^*_{\beta}} g(M)-\kappa^{-p^*_{\beta}} g(\kappa)\right)
\\[4mm]\phantom{--------------\displaystyle-\int_{\kappa}^M  g(s)d\omega(s) }\displaystyle
+\frac{c_{14}}{p^*_{\beta}+1}\int_{\kappa}^M
s^{-1-p^*_{\beta}} g(s)ds.
\end{array}$$
Since  $\lim_{M\to\infty}M^{-p^*_{\beta}} g(M)=0$ by  $(\ref{1.4})$ and \cite[Lemma 4.1]{CV2} and
$\omega(s)\leq c_{14}s^{-p^*_{\beta}}$, we derive
 $ g(\kappa)\omega(\kappa)\leq c_{14}\kappa^{-p^*_{\beta}} g(\kappa)$ and then
$$-\int_{\kappa}^\infty  g(s)d\omega(s)\leq \frac{c_{14}}{p^*_{\beta}+1}\int_{\kappa}^\infty s^{-1-p^*_{\beta}} g(s)ds.
$$
The above quantity on the right-hand side tends to $0$
when $\kappa\to\infty$. The conclusion follows: for any
$\epsilon>0$ there exists $\kappa>0$ such that
$$\frac{c_{14}}{p^*_{\beta}+1}\int_{\kappa}^\infty s^{-1-p^*_{\beta}} g(s)ds\leq \frac{\epsilon}{2}
$$
and there exists $\delta>0$ such that
$$\int_E t^{\beta}dxdt\leq \delta\Longrightarrow  g(\kappa)\int_E t^{\beta}dxdt\leq\frac{\epsilon}{2}.
$$
This means that $\{h_{n_k}\circ u_{\nu_{n_k }}\}$ is uniformly integrable in
$L^1(Q_T)$ and by Vitali convergence theorem  $h_{n_k}\circ u_{\nu_{n_k }}\to h\circ U$ in
$L^1(Q_T)$ . Letting
$n_k\to\infty$ in the identity
$$u_{\nu_{n_k }}(t,.)=\mathbb H_{\alpha}[\nu_{n_k}](t,.)-\myint{0}{t}\mathbb H_{\alpha}[h\circ u_{\nu_{n_k }}(s,.)](t-s,.)ds
$$
for some $t>0$ such that $u_{\nu_{n_k }}(t,.)\to U(t,.)$ a.e. in $\R^N$ yields
$$U(t,.)=\mathbb H_{\alpha}[\nu](t,.)-\myint{0}{t}\mathbb H_{\alpha}[h\circ U(s,.)](t-s,.)ds.
$$
This is valid for almost all $t>0$ and implies that $U\in C([0,T];L^1(\R^N))$, up to a modification on a set of $t>0$ with zero measure. Moreover
$$\begin{array}{lll}
\displaystyle\int_{Q_T}\left(u_{\nu_{n_k }}(-\partial_t\xi+(-\Delta)^{\alpha}\xi)+\xi h\circ u_{\nu_{n_k }}\right) dxdt
\\\phantom{-----------}\displaystyle=\int_{\R^N}\xi(0,x)\nu_{n_k } dx-\int_{\R^N}u_{\nu_{n_k }}(T,x)\xi(T,x) dx.
\end{array}$$
where $\xi\in \mathbb Y_{\alpha,T}$ is arbitrary. Thus, using the continuity of $t\mapsto U(t,.)$ in $L^1(\R^N)$, we derive
$$\begin{array}{lll}
\displaystyle\int_{Q_T}\left(U(-\partial_t\xi+(-\Delta)^{\alpha}\xi)+\xi h\circ U\right) dxdt
\\\phantom{-----------}\displaystyle=\int_{\R^N}\xi(0,x)d\nu(x)-\int_{\R^N}U(T,x)\xi(T,x) dx.
\end{array}$$
From this we infer that $U$ is a weak
solution of $(\ref{he 1.1})$.\smallskip

\noindent{\it  Existence for general $\nu$. } For $\nu\in\mathfrak{M}^b(\R^N)$, a sequence $\{\nu_n\}$ in $C^2_0(\R^N)$
converge to $\nu$ in the weak sense of bounded measures. Because of the monotonicity of $h(t,\cdot)$,
$$-\mathbb{H}_{\alpha}[|\nu_n|]\le u_{-|\nu_n|}\le u_{\nu_n}\le  u_{|\nu_n|}\le \mathbb{H}_{\alpha}[|\nu_n|].$$
Then by above analysis, the sequence $\{h\circ u_{-|\nu_n|})\}$ and $\{h\circ u_{|\nu_n|})\}$ are relatively compact in $L^1(Q^B_T)$ for any
$T>0$ and ball $B$ and (\ref{L8}) holds for $\{u_{\nu_n}\}$. Therefore  $\{u_{\nu_n}\}$ is relatively locally compact in $L^1(Q^B_T)$ and there exist some
subsequence $\{u_{\nu_{n_k}}\}$ and $U\in L^1(Q_T)$ such that
$$u_{\nu_{n_k}}\to U\Longrightarrow h\circ u_{\nu_{n_k}}\to  h\circ U\quad {\rm as}\;\, k\to\infty\quad {\rm a.e.\;  in}\quad Q_T.$$
As in the previous case it implies that $U$ is a weak solution of $(\ref{he 1.1})$ and also an integral solution.\smallskip

\noindent{\it   Uniqueness.} Let $u_1,u_2$ be two weak solutions of
(\ref{he 1.1}) with the same initial $\nu$ and $w=u_1-u_2$. Then
$$\partial_tw+(-\Delta)^\alpha w=h\circ u_2-h\circ u_1\quad{\rm in}\ \ Q_T.$$
Since $h\circ u_2-h\circ u_1\in L^1(Q_T)$, then by (\ref{sign}),  for $\xi\in\mathbb{Y}_{\alpha,T}$, $\xi\ge0$, we have that
\begin{eqnarray*}
&&\displaystyle\int_{Q_T} |w|[-\partial_t\xi+ (-\Delta)^\alpha \xi]
dxdt+\int_{\R^N}|w(T,x)|\xi(T,x)dxdt\\&&\quad \quad\quad \quad \quad\quad \displaystyle+\int_{Q_T}(h\circ u_2-h\circ u_1){\rm sign}(w)\xi dxdt\le0.
\end{eqnarray*}
This implies $w=0$ by monotonicity.\smallskip

\noindent Statements (i) and (ii) and inequality (\ref{12-09-2}) follows from the fact that the same relation holds for $u_{\nu_n}$ by
Lemma \ref{lm 2}.
\smallskip

\noindent Stability is proved by the same approach that existence. If $\{\nu_n\}$ converges to $\nu$ in the weak sense of measures, then $\norm{\nu_n}_{\mathfrak M^b}$ is bounded independently of $n$. Since the distribution function of $h\circ u_{\nu_n}$ depends only on the supremum of $\norm{\nu_n}_{\mathfrak M^b}$, this set of functions is uniformly integrable in $Q_T$. This, combined with local compactness of the set $\{u_{\nu_n}\}$ in $L^1(Q_T)$, implies the convergence of a subsequence
$(u_{\nu_{n_k}},h\circ u_{\nu_{n_k}})$ to $(u_\nu,h\circ u_\nu)$ where $u_\nu$ is the solution of (\ref{he 1.1}). Because of uniqueness, all converging subsequences have the same limit, which imply the convergence of the whole sequence and  stability.
 \hfill$\Box$\medskip

\setcounter{equation}{0}
\section{Dirac mass as initial data}

In this section, we study the properties of solutions to (\ref{he 1.1}) when $h(t,r)=t^\beta r^p$
with $\beta>-1$ and $0<p<p^*_\beta$ and the initial data is $\nu=k\delta_0$ with $k>0$.

\begin{proposition}\label{pr 4.2}
Assume $0<p<p^*_\beta$ and that  $u_k$ is the solution of (\ref{2.1}),
then there exists $c_{15}>0$ such that
\begin{equation}\label{13-009-0}
\lim_{t\to0^+} t^{\frac N{2\alpha}}u_k(t,0)=c_{15}k.
\end{equation}

\end{proposition}
{\bf Proof.}  By (\ref{12-09-2}) it follows that
\begin{equation}\label{12-09-4}
  u_k(t,0)\le k\mathbb{H}_\alpha[\delta_0](t,0)=k\Gamma_\alpha(t,0),\qquad t>0.
\end{equation}
We claim that
there exists  $c_{16}>0$ independent of $k$ such that
\begin{equation}\label{26-09-1-1}
u_k(t,0)\ge k\Gamma_\alpha(t,0)-c_{16}k^p t^{-\frac N{2\alpha}p+1+\beta}, \quad t\in(0,1/2).
\end{equation}
Indeed, from (\ref{12-09-2}), it is infered that
$$u_k(t,0)\ge k\Gamma_\alpha(t,0) -k^pW(t,0),\qquad t\in(0,1/2),$$
where
$$
W(t,x)=\myint{0}{t}\mathbb{H}_\alpha[s^\beta(\mathbb{H}^p_\alpha[\delta_0]](t-s,x)ds,\qquad (t,x)\in Q_\infty.
$$
For $t\in(0,1/4)$, there exists $c_{17},c_{18}>0$ such that
$$
\begin{array} {ll}
\displaystyle
   W(t,0)\le c_{17}\int_0^t\int_{\R^N}\frac{(t-s)^{-\frac{N}{2\alpha}}s^\beta}{1+((t-s)^{-\frac{1}{2\alpha}}|y|)^{N+2\alpha}}
  \left(\frac{s^{-\frac{N}{2\alpha}}}{1+(s^{-\frac{1}{2\alpha}}|y|)^{N+2\alpha}}\right)^pdyds
   \\[4mm]\displaystyle\phantom{---\ }
    \le c_{17}\int_0^t\int_{\R^N}\frac{s^{\beta-\frac N{2\alpha}p}dzds}{\left(1+\left((\frac{t-s}s)^{\frac 1{2\alpha}}|z|\right)^{(N+2\alpha)p}\right)\left(1+|z|^{N+2\alpha}\right)}
   \\[6mm]\displaystyle\phantom{---\ }
   \le c_{17}t^{\beta+1-\frac{Np}{2\alpha}}\int_0^1\int_{\R^N}
\frac{d\tau dZ}{\left(1+\left(\frac{1-\tau}{\tau}\right)^{\frac{(N+2\alpha)p}{2\alpha}}|Z|^{(N+2\alpha)p}\right)\left(1+|Z|^{N+2\alpha}\right)}

    \\[4mm]\displaystyle\phantom{---\ }
    \le c_{18}t^{\beta+1-\frac{Np}{2\alpha}}.
\end{array}
$$
Combining  (\ref{17-10-0}) and  $-\frac{N}{2\alpha}p+1+\beta>-\frac N{2\alpha}$, we obtain that
$$\lim_{t\to0^+}t^{\frac N{2\alpha}}W(t,0)=0. $$
Therefore, (\ref{13-009-0}) holds.\hfill$\Box$
\smallskip

 In what follows we consider the limit of the solution $\{u_k\}$ of (\ref{2.1}) as $k\to\infty$ for $p\in(0,1]$.

\begin{proposition}\label{prop4.2}
Assume $0<p\leq 1$ and  that $u_k$ is the solution of (\ref{2.1}), then
$$\lim_{k\to\infty}u_k=\infty\quad{\rm in}\quad   Q_\infty,$$
locally uniformly in $Q_\infty$.
\end{proposition}
{\bf Proof.}  We observe that $\mathbb{H}_{\alpha}[\delta_0]$ and $ \mathbb{H}_{\alpha}[t^\beta(\mathbb{H}_{\alpha}[\delta_0])^p]$ are positive in $(0,\infty)\times\R^N$.
By (\ref{12-09-2}), for $p\in(0,1)$ and $(t,x)\in  (0,\infty)\times \R^N$, we have that
\begin{eqnarray*}\displaystyle
u_k&\ge& k\mathbb{H}_{\alpha}[\delta_0]- k^pW\Longrightarrow \lim_{k\to\infty}u_k=\infty.
\end{eqnarray*}
For $p=1$, it is obvious that $u_k=ku_1$ and $u_1>0$ in $(0,\infty)\times\R^N$, then
$$\lim_{k\to\infty}u_k=\infty\quad {\rm{in}}\quad  Q_\infty. $$
The proof is complete. \hfill$\Box$ \smallskip

Now we deal with the range $p\in(1,p^*_\beta)$.

\begin{lemma}\label{lm 4.1}
Assume $1<p<p^*_\beta$ and that $u_k$ is the solution of (\ref{2.1}).
Then for any $k>0$,
\begin{equation}\label{13-09-0}
  0\le u_k\le U_p\quad{\rm in}\ \ Q_\infty,
\end{equation}
where $U_p$ is given by (\ref{20-10-0}).
\end{lemma}
{\bf Proof.} Let $\{f_{n,k}\}$ be a sequence of nonnegative functions in $C^1_c(\R^N)$
which converges to $k\delta_0$ in the weak sense of measures as $n\to\infty$. We denote by $u_{n,k}$ the corresponding solution
of (\ref{he 2.1}) with initial data by $f_{n,k}$.\smallskip

\noindent{\it We claim that}
\begin{equation}\label{4.1}
u_{n,k}\le U_p \quad {\rm in}\quad Q_\infty,
\end{equation}
where, we recall it, $U_p$ is the maximal solution of the ODE $y'+t^\beta y^p=0$ on $\R_+$. Indeed this implies (\ref{13-09-0}).
\smallskip

\noindent\emph{Step 1. We claim that }
\begin{equation}\label{4.0}
\lim_{|x|\to\infty}u_{n,k}(t,x)=0,\qquad\forall t>0.
\end{equation}
From  \cite{CKS,CT}, there exists $c_8>0$ such that for any $x,y\in\R^N$ and $t\in(0,\infty)$,
$$0<\Gamma_\alpha(t,x-y)\le\frac{c_8t^{-\frac N{2\alpha}}}{1+(|x-y|t^{-\frac1{2\alpha}})^{N+2\alpha}}.$$
Then for $|x|>1$,
\begin{eqnarray*}
0\le\mathbb{H}_\alpha[f_{n,k}](t,x) &\le &c_8t^{-\frac N{2\alpha}}\int_{\R^N}\frac{f_{n,k}(y)}{1+(|x-y|t^{-\frac1{2\alpha}})^{N+2\alpha}}dy
  \\&=& c_8\int_{\R^N}\frac{f_{n,k}(x-zt^{\frac1{2\alpha}})}{1+|z|^{N+2\alpha}}dz
  \\&=&  c_8\left(\int_{\R^N\setminus B_R}\frac{f_{n,k}(x-zt^{\frac1{2\alpha}})}{1+|z|^{N+2\alpha}}dz+  \int_{B_R}\frac{f_{n,k}(x-zt^{\frac1{2\alpha}})}{1+|z|^{N+2\alpha}}dz\right),
\end{eqnarray*}
where $R=\frac{1}2|x|t^{-\frac1{2\alpha}}$  and $B_R=\{z\in\R^N:|z|<R\}$.
It is obvious that
$$|x-zt^{\frac1{2\alpha}}|\ge|x|-|z|t^{\frac1{2\alpha}}\geq  |x|/2\quad \ {\rm for\ all}\ z\in B_R. $$
Then
\begin{eqnarray*}
\int_{B_R}\frac{f_{n,k}(x-zt^{\frac1{2\alpha}})}{1+|z|^{N+2\alpha}}dz&\le& \sup_{|y|\ge \frac{|x|}2}f_{n,k}(y)\int_{B_R}\frac{1}{1+|z|^{N+2\alpha}}dz
\\ &\le& \sup_{|y|\ge \frac{|x|}2}f_{n,k}(y)\int_{\R^N}\frac1{1+|z|^{N+2\alpha}}dz
\\&=& c_{16} \sup_{|y|\ge \frac{|x|}2}f_{n,k}(y)
\end{eqnarray*}
and
\begin{eqnarray*}
\int_{\R^N\setminus B_R}\frac{f_{n,k}(x-zt^{\frac1{2\alpha}})}{1+|z|^{N+2\alpha}}dz
 \le  \int_{\R^N\setminus B_R}\frac{\|f_{n,k}\|_{L^\infty(\R^N)}}{1+|z|^{N+2\alpha}}dz\leq c_{18}R^{-2\alpha}=\frac {c_{18}t}{ |x|^{2\alpha}},
\end{eqnarray*}
for some  $c_{18}>0$ independent of $x,t$ and $R$.
Since $f_{n,k}\in C_0^1(\R^N)$, we have that
$$\lim_{|x|\to\infty}\sup_{|y|\ge \frac{|x|}2}f_{n,k}(y)=0$$
 and then for any $t>0$,
$0\le u_{n,k}(t,x)\le \mathbb{H}_\alpha[f_{n,k}](t,x) \to 0$ as $|x|\to\infty$.

\smallskip

\noindent\emph{Step 2. We claim that (\ref{4.1}) holds.} By contradiction, if (\ref{4.1}) is not verified,  there exists $(t_0,x_0)\in (0,\infty)\times\R^N$ such that
$$( U_p-u_{n,k})(t_0,x_0)=\min_{(t,x)\in (0,\infty)\times\R^N}( U_p-u_{n,k})(t,x)<0,$$
since $ U_p(t)>0=\lim_{|x|\to\infty}u_{n,k}(t,x)$ for any $t\in (0,\infty)$, $U_p(0)=\infty>f_{n,k}(x)=u_{n,k}(0,x)$ for $x\in\R^N$
and $\lim_{t\to\infty}U_p(t)=\lim_{t\to\infty}u_{n,k}(t,x)=0$ for $x\in\R^N$.
Then  $\partial_t( U_p-u_{n,k})(t_0,x_0)=0$. Moreover, 
$$\begin{array}{lll}
( U_p-u_{n,k})(t_0,x_0)=\min\{U_p(t_0)-u_{n,k}(t_0,x):x\in\R^N\}\\[2mm]
\phantom{( U_p-u_{n,k})(t_0,x_0)}
=U_p(t_0)-\max\{u_{n,k}(t_0,x):x\in\R^N\}
\end{array}$$
and
$$u_{n,k}(t_0,x_0)=\max\{u_{n,k}(t_0,x):x\in\R^N\}\Longrightarrow (-\Delta)^\alpha  u_{n,k} (t_0,x_0)\ge 0.$$
Then
\begin{eqnarray*}
0=\partial_t( U_p-u_{n,k})(t_0,x_0)-(-\Delta)^\alpha u_{n,k}(t_0,x_0)+ t_0^\beta U_p^p(t_0)-t_0^\beta u_{n,k}^p(t_0,x_0)<0,
\end{eqnarray*}
which is impossible. Thus (\ref{4.1}) holds. \hfill $\Box$\medskip


\begin{proposition}\label{re 12-09}
(i) Assume $0<p<p^*_\beta$ and  that $u_k$ is the solution of (\ref{2.1}).
Then $u_k$ is a classical solution of (\ref{he 2.1}).\smallskip

\noindent (ii) Assume $1<p<p^*_\beta$ and that $u_\infty$ is defined  by (\ref{13-09-080}).
Then $u_\infty$ is a classical solution of (\ref{he 2.1}).
\end{proposition}
{\bf Proof.} $(i)$
Since $u_k\le k\mathbb{H}_\alpha[\delta_0]$, it is infered that $u_k$ is  bounded in $(\epsilon,\infty)\times\R^N$ for $\epsilon>0$.
Let $\{g_{n,k}\}$ be a sequence of nonnegative functions in $C^1_0(\R^N)$
which converges to $k\delta_0$ as $n\to\infty$  and $u_{n,k}$ the corresponding solution
of (\ref{he 2.1}) with initial data $g_{n,k}$. Then $\mathbb{H}_\alpha[g_{n,k}]\to k\mathbb{H}_\alpha[\delta_0]$ as $n\to\infty$ uniformly in $[\epsilon,\infty)\times\R^N$ for any $\epsilon>0$ and by the Comparison Principle, there exists $c_{19}>1$ such that
$$0\le u_{n,k}(t,x)\le k\mathbb{H}_\alpha[g_{n,k}]\le c_{19}k \mathbb{H}_\alpha[\delta_0]\quad {\rm in}\ \ [\epsilon,\infty)\times\R^N,$$
and there exists $\sigma\in(0,1)$ such that
$\{u_{n,k}\}$ are uniformly bounded with respect to $n$ in $C^{\frac\sigma{2\alpha},\sigma}_{t,x}((\epsilon,\infty)\times\R^N)$ with $\epsilon>0$.
Therefore, by the Arzela-Ascoli theorem, $u_{n,k}$ converges to
$u_k$ in $C^{\frac{\sigma'}{2\alpha},\sigma'}_{t,x}((\epsilon,\infty)\times\R^N)$ with $\sigma'\in(0,\sigma)$ and then $u_k$  is a viscosity solution of (\ref{he 2.1}) in $(\epsilon,\infty)\times\R^N$. By  estimate $(A.1)$ in \cite{CF},  $u_k$ is in $C^{1+\sigma',2\alpha+\sigma'}_{t,x}((\epsilon,\infty)\times\R^N)$ and $u_k$ is a classical solution of
(\ref{he 2.1}) in $(\epsilon,\infty)\times\R^N$.\smallskip

 \noindent (ii)  The proof is the same as part $(i)$, just replacing  $u_k\le k\mathbb{H}_\alpha[\delta_0]$ by $u_\infty\le U_p$. \hfill$\Box$

\setcounter{equation}{0}
\section{Self-similar and very singular solutions}

By Theorem \ref{teo 1} and (\ref{13-09-0}), we see that $\{u_k\}$ is an increasing sequence of nonnegative functions bounded from above by $U_p$.
Then for $p\in(1,p_\beta^*)$, there exists $u_\infty=\lim_{k\to\infty}u_k$, which is a classical solution of
(\ref{he 2.1})
by Proposition \ref{re 12-09} $(ii)$ and satisfies
\begin{equation}\label{02-10-3}
u_\infty\le U_p\quad{\rm in}\quad Q_\infty.
\end{equation}

\begin{proposition}\label{lm 13-09-1}
Assume $1<p<p_\beta^*$, then $u_\infty$ is a self-similar solution of (\ref{he 2.1}).
\end{proposition}
{\bf Proof.} For $\lambda>0$, we set
$$T_{\lambda} [u](t,x)=\lambda^{\frac{2\alpha(1+\beta)}{p-1}}u(\lambda^{2\alpha}t,\lambda x),\qquad (t,x)\in Q_\infty. $$
It is straightforward to verify that $T_{\lambda} [u_k]$ is the solution of
\begin{equation}\label{scal}\begin{array}{lll}
\partial_t u + (-\Delta)^\alpha  u +t^\beta u^p=0\qquad&\text{in } Q_\infty,\\
\phantom{\partial_t  +., -\Delta^\alpha   t^\beta u^p}
u(0,.)=\lambda^{\frac{2\alpha(1+\beta)}{p-1}-N}k\delta_0\qquad&\text{in } \R^N.
\end{array}
\end{equation}
Because of uniqueness,  $T_{\lambda} [u_k]=u_{k\lambda^{\frac{2\alpha(1+\beta)}{p-1}-N}}$.
Letting $k\to\infty$ and using the continuity of $u\mapsto T_\lambda[u]$, we have that
$$
\lim_{k\to\infty}T_{\lambda} [u_k]=T_{\lambda} [u_\infty]=u_\infty,
$$
which implies that $u_\infty$ is a self-similar solution (\ref{he 2.1}).
\hfill$\Box$ \medskip

Let us denote
$$U_\infty(z)=u_\infty(1,z),\qquad z\in\R^N, $$
then
$U_\infty$ is a classical solution of (\ref{13-09-07}).
It is clear  that the constant $(\frac{1+\beta}{p-1})^{\frac1{p-1}}$ is a constant  positive solution of the self-similar equation (\ref{13-09-07}).
We observe that $N<\frac{2\alpha(1+\beta)}{p-1}<N+2\alpha$ when $1+\frac{2\alpha(1+\beta)}{N+2\alpha}<p<1+\frac{2\alpha(1+\beta)}{N}$.

We prove below this fundamental result that   $u_\infty$ is the minimal self similar solution.

\begin{proposition}\label{lm 5.1}
Assume that $1<p<1+\frac{2\alpha(1+\beta)}{N}$ and $\tilde u$  is a positive self-similar solution of (\ref{5.1}). Then $u_\infty\leq \tilde u$.
\end{proposition}
{\bf Proof.} For any $r>0$, we have that
\begin{eqnarray*}
 \int_{B_r(0)}\tilde u(t,x)dx &=& t^{-\frac{1+\beta}{p-1}} \int_{B_r(0)}\tilde u(1,t^{-\frac1{2\alpha}}x)dx
 \\&=& t^{-\frac{1+\beta}{p-1}+\frac N{2\alpha}}\int_{B_{t^{-\frac1{2\alpha}}r}(0)}\tilde u(1,z)dz
 \\&\ge& t^{-\frac{1+\beta}{p-1}+\frac N{2\alpha}}\int_{B_{1}(0)}\tilde u(1,z)dz
 \\&\to&+\infty\quad{\rm as}\ t\to0^+,
\end{eqnarray*}
where last inequality holds for $t\in (0,r^{2\alpha}]$. Let $\{\epsilon_n\}$ be a sequence positive decreasing numbers converging to 0 as $n\to\infty$.
 For $\epsilon_n$ and $k>0$, there exists $t_{n,k}>0$ such   that
$$\int_{B_{\epsilon_n}(0)}\tilde u(t_{n,k},x)dx=k.$$
We observe that for any fixed $k$, $t_{n,k}\to0$ as $n\to\infty$ since $\lim_{n\to\infty}\epsilon_n=0$.
Let $\eta_0:\R^N\to[0,1]$ be a $C^2$ function such that supp$\,\eta_0\subset \bar B_2(0)$, $\eta_0=1$ in $B_1(0)$ and $\eta_n(x)=\eta_0(\epsilon_n^{-1}x)$ for
$x\in\R^N$.
Choosing  $\{f_{n,k}\}$ be a sequence of $C^2$ functions such that
 $$0\le f_{n,k}(x)\le\eta_n(x)\tilde u(t_{n,k},x),\qquad\forall  x\in\R^N$$
and
\begin{eqnarray*}
f_{n,k}\to k\delta_0\qquad {\rm as}\quad  n\to\infty.
\end{eqnarray*}
Let $u_{n,k}$ be the solution of (\ref{he 1.1}) with initial data $f_{n,k}$,
then $$u_{n,k}(t,x)\le u(t_{n,k}+t,x),\qquad\forall  (t,x)\in Q_\infty$$
 and by uniqueness of $u_k$,
$\lim_{n\to\infty}u_{n,k}=u_k$, where
$u_k$ is the solution of (\ref{he 1.1}) with initial data $k\delta_0$.
Then for any $k$, we have
$u_k\le \tilde u$ in $Q_\infty$,
which implies that
$$ u_\infty\le\tilde u\quad {\rm in}\quad Q_\infty.$$
 \hfill$\Box$

\subsection{The case $1+\frac{2\alpha(1+\beta)}{N+2\alpha}<p<1+\frac{2\alpha(1+\beta)}{N}$ }


We define the function $ w_\lambda$ by
\begin{equation}\label{17-09-0}
 w_\lambda(t,x)=\lambda t^{-\frac{1+\beta}{p-1}}w(t^{-\frac1{2\alpha}}|x|),\qquad  (t,x)\in Q_\infty,
 \end{equation}
where $w(s)=\frac{\ln(e+s^2)}{1+s^{N+2\alpha}}$.

\begin{lemma}\label{lm 14-09-0}
Assume $1+\frac{2\alpha(1+\beta)}{N+2\alpha}<p<1+\frac{2\alpha(1+\beta)}{N}$, then there exists $\Lambda_0>0$ such that for $\lambda\ge\Lambda_0$,
\begin{equation}\label{14-09-0}
\partial_tw_\lambda(t,x)+(-\Delta)^\alpha w_\lambda(t,x)+t^\beta w_\lambda^p(t,x)\ge0,\quad \forall(t,x)\in Q_\infty.
\end{equation}

\end{lemma}
{\bf Proof.}  By direct computation, we have
$$
\partial_tw_\lambda(t,x) = -\frac{\lambda(1+\beta)}{p-1}t^{-\frac{1+\beta}{p-1}-1}w(t^{-\frac1{2\alpha}}|x|)
-\frac{\lambda}{2\alpha}t^{-\frac{1+\beta}{p-1}-\frac1{2\alpha}-1}|x|w'(t^{-\frac1{2\alpha}}|x|)
$$
and
$$(-\Delta)^\alpha w_\lambda(t,x)=\lambda t^{-\frac{1+\beta}{p-1}-1}(-\Delta)^\alpha w(t^{-\frac1{2\alpha}}|x|),$$
which implies
\begin{equation}\label{14-09-00}\begin{array}{lll}
\partial_t w_\lambda(t,x)+(-\Delta)^\alpha w_\lambda(t,x)+t^\beta w^p_\lambda(t,x) \\[2mm]
\phantom{---}
= \lambda t^{-\frac{1+\beta}{p-1}-1}\left[(-\Delta)^\alpha w(s)-\myfrac1{2\alpha}w'(s)s-\myfrac{1+\beta}{p-1}w(s)+\lambda^{p-1}w^p(s)\right] ,
\end{array}\end{equation}
where $s=|z|$ with $z=t^{-\frac1{2\alpha}}x$. Next, for $s>0$, we have
$$
-\frac1{2\alpha}w'(s)s-\frac{1+\beta}{p-1}w(s) = \left[\frac{N+2\alpha}{2\alpha}\frac{s^{N+2\alpha}}{1+s^{N+2\alpha}}-\frac{1+\beta}{p-1}
-\frac{s^2(e+s^2)^{-1}}{\alpha\ln(e+s^2)}\right]w(s).
$$
Since $\frac{N+2\alpha}{2\alpha}>\frac{1+\beta}{p-1}$, $\lim_{s\to\infty}\frac{s^{N+2\alpha}}{1+s^{N+2\alpha}}=1$ and $\lim_{s\to\infty}\frac1{\ln(e+s^2)}=0$,
there exists $R_0>0$ and $\sigma_0>0$ such that
\begin{equation}\label{14-09-5}
-\frac1{2\alpha}w'(s)s-\frac{1+\beta}{p-1}w(s)\ge \sigma_0 w(s),\qquad \forall s\ge R_0.
\end{equation}
  For $|z|>2$, and using the definition of the fractional Laplacian, we have
\begin{equation}\label{15-09-0}\begin{array}{lll}\displaystyle
-(-\Delta)^\alpha w(|z|) =\frac{1}{2}\int_{\R^N}\left(\frac{\ln(e+|z+\tilde y|^2)}{1+|z+\tilde y|^{N+2\alpha}}
+\frac{\ln(e+|z-\tilde y|^2)}{1+|z-\tilde y|^{N+2\alpha}}-\frac{2\ln(e+|z|^2)}{1+|z|^{N+2\alpha}}\right)
\myfrac{d \tilde y}{|\tilde y|^{N+2\alpha}}
\\[4mm]\phantom{-(-\Delta)^\alpha w(|z|)}\displaystyle
=\frac{w(|z|)}{2|z|^{2\alpha}}\int_{\R^N}\frac{I_z(y)}{|y|^{N+2\alpha}}dy,
\end{array}\end{equation}
where
$$\begin{array}{lll}\displaystyle
 I_z(y) = \frac{1+|z|^{N+2\alpha}}{1+|z|^{N+2\alpha}|e_z+y|^{N+2\alpha}}\frac{\ln(e+|z|^2|e_z+y|^2)}{\ln(e+|z|^2)} \\[2mm]
\phantom{ I_z(y)----} \displaystyle
+\frac{1+|z|^{N+2\alpha}}{1+|z|^{N+2\alpha}|e_z-y|^{N+2\alpha}} \frac{\ln(e+|z|^2|e_z-y|^2)}{\ln(e+|z|^2)}-2
\end{array}$$
 and $e_z=\frac{z}{|z|}$.\smallskip

\noindent \emph{We claim that there exists $c_{20}>0$ such that
 \begin{equation}\label{15-09-1}
   \int_{B_{\frac12}(-e_z)\cup B_{\frac12}(e_z)}\frac{I_z(y)}{|y|^{N+2\alpha}}dy\le \frac{ c_{20}}{w(|z|)|z|^N}.
 \end{equation}}
In fact, for $y\in B_{\frac12}(-e_z)$, there exists $c_{21}>0$ such that
$$\frac{1+|z|^{N+2\alpha}}{1+|z|^{N+2\alpha}|e_z-y|^{N+2\alpha}} \frac{\ln(e+|z|^2|e_z-y|^2)}{\ln(e+|z|^2)}\le c_{21}$$
and  then
\begin{eqnarray*}
  \int_{B_{\frac12}(-e_z)}\frac{I_z(y)}{|y|^{N+2\alpha}}dy &\le&
  \omega_N\int_0^{\frac12}\frac{1+|z|^{N+2\alpha}}{1+(|z| r)^{N+2\alpha}}\frac{\ln(e+|z|^2 r^2)}{\ln(e+|z|^2)}r^{N-1}dr+c_{22}
    \\&\le&\frac{\omega_N}{w(|z|)|z|^N}\int_0^{\infty}\frac{t^{N-1}\ln(e+t^2)}{1+ t^{N+2\alpha}}dt+c_{22}
    \\&\le&\frac{c_{23}}{w(|z|)|z|^N},
\end{eqnarray*}
where $c_{22},c_{23}>0$ and the last inequality holds since $w(|z|)|z|^N\to0$ as $|z|\to\infty$. Thus,

$$\int_{B_{\frac12}(e_z)}\frac{I_z(y)}{|y|^{N+2\alpha}}dy=\int_{B_{\frac12}(-e_z)}\frac{I_z(y)}{|y|^{N+2\alpha}}dy\le \frac{c_{23}}{w(|z|)|z|^N}.$$\smallskip

\noindent \emph{ We claim that there exists $c_{24}>0$ such that
 \begin{equation}\label{15-09-2}
   \int_{B_{\frac12}(0)}\frac{I_z(y)}{|y|^{N+2\alpha}}dy\le c_{24}.
 \end{equation}}
Indeed, since the function $I_z$ is $C^2$ in $\bar B_{\frac12}(0)$, $I_z(0)=0$ and   $I_z(y)=I_z(-y),$
then $\nabla I_z(0)=0$ and there exists $c_{34}>0$ such that
$$|D^2 I_z(y)|\le c_{25}\qquad \forall y\in B_{\frac12}(0). $$
Then we have
$$I_z(y)\le c_{25}|y|^2\qquad \forall y\in B_{\frac12}(0),$$
which implies
$$ \int_{B_{\frac12}(0)}\frac{I_z(y)}{|y|^{N+2\alpha}}dy\le c_{25}\int_{B_{\frac12}(0)}\frac{|y|^{2}}{|y|^{N+2\alpha}}dy\le c_{24}.
$$\smallskip

\noindent \emph{ We claim that there exists $c_{26}>0$ such that
 \begin{equation}\label{15-09-3}
   \int_{A}\frac{I_z(y)}{|y|^{N+2\alpha}}dy\le c_{26},
 \end{equation}
where $A=\R^N\setminus (B_{\frac12}(0)\cup B_{\frac12}(e_z)\cup B_{\frac12}(-e_z))$.}
In fact, for $y\in A$, we observe that there exists $c_{27}>0$ such that
$I_z(y)\le c_{27}$ and
$$
   \int_{A}\frac{I_z(y)}{|y|^{N+2\alpha}}dy\le \int_{\R^N\setminus B_{\frac12}(0)}\frac{c_{27}}{|y|^{N+2\alpha}}\le c_{28},
$$
for some $c_{28}>0$.
Therefore, by (\ref{14-09-00})-(\ref{15-09-3}), there exists $c_{29}>0$ such that
\begin{equation}\label{15-09-4}
 (-\Delta)^\alpha w(|z|)\ge -\frac{c_{29}}{1+|z|^{N+2\alpha}},\qquad |z|\ge 2.
\end{equation}

 \noindent By (\ref{14-09-5}) and (\ref{15-09-4}), there exists  $R_1\ge R_0+2$ such that for $|z|>R_1$,
$$
\begin{array} {ll}
\displaystyle
   (-\Delta)^\alpha w(|z|)-\frac1{2\alpha}w'(|z|)|z|-\frac{1+\beta}{p-1}w(|z|)
   \ge\sigma_0w(|z|)-\frac{c_{29}}{1+|z|^{N+2\alpha}}
   \\[4mm]\displaystyle\phantom{(-\Delta)^\alpha w(|z|)-\frac1{2\alpha}w'(|z|)|z|-\frac{1+\beta}{p-1}w(|z|)}
   =w(|z|)\left(\sigma_0-\frac{c_{29}}{\ln(e+|z|^2)}\right)
     \\[4mm]\displaystyle\phantom{(-\Delta)^\alpha w(|z|)-\frac1{2\alpha}w'(|z|)|z|-\frac{1+\beta}{p-1}w(|z|)}
   \ge0.
\end{array}
$$
When $|z|\le R_1$, it is clear that there exists $c_{30}>0$ such that
$$(-\Delta)^\alpha w(|z|)-\frac1{2\alpha}w'(|z|)|z|-\frac{1+\beta}{p-1}w(|z|)\ge -c_{30}.$$
Then there exists $\Lambda_0>0$ such that for $\lambda\ge\Lambda_0$,
\begin{equation}\label{15-09-5}
  (-\Delta)^\alpha w(|z|)-\frac1{2\alpha}w'(|z|)|z|-\frac{1+\beta}{p-1}w(|z|)+\lambda^{p-1}w^p(|z|)\ge0,\qquad \forall z\in\R^N,
\end{equation}
which, together with (\ref{14-09-00}), implies that (\ref{14-09-0}) holds.
\hfill$\Box$\medskip

Next we prove that $u_\infty$ is not a trivial flat solution when $1+\frac{2\alpha(1+\beta)}{N+2\alpha}<p<p^*_\beta$.

\begin{lemma}\label{lm 13-09-2}
Assume $1+\frac{2\alpha(1+\beta)}{N+2\alpha}<p<1+\frac{2\alpha(1+\beta)}{N}$, that $w_{\Lambda_0}$ is given in (\ref{17-09-0})
and $u_\infty$ is given in (\ref{13-09-080}).
Then
\begin{equation}\label{13-09-9}
u_\infty(t,x)\le w_{\Lambda_0}(t,x)\qquad \forall (t,x)\in Q_\infty.
\end{equation}
Moreover,
\begin{equation}\label{13-09-90}
\lim_{t\to 0}u_\infty(t,\cdot)=0\quad \text {uniformly on }\; B_\epsilon^c,\quad\forall \epsilon >0.
\end{equation}

\end{lemma}
{\bf Proof.}
 Let us denote
 $$f_0(r)=\frac{k_0\ln(e+r^2)}{1+r^{N+2\alpha}},\quad\forall\;  r\ge0\quad {\rm and}\quad f_{n,k}(x)=kn^Nf_0(n|x|),\quad\forall   x\in\R^N,$$
 where  $$k_0=\left[\omega_N\int_0^\infty\frac{\ln(e+r^2)}{1+r^{N+2\alpha}}r^{N-1}dr\right]^{-1}.$$
Then for any $\eta\in C_c(\R^N)$, we have that
\begin{eqnarray*}
 \lim_{n\to\infty}\int_{\R^N} f_{n,k}\eta dx = k\lim_{n\to\infty}\int_{\R^N}f_0(|x|)\eta\left(\frac{x}{n}\right)dx = k\eta(0).
\end{eqnarray*}
Let $t_n=n^{-2\alpha}$ and then
\begin{eqnarray*}
 w_{\Lambda_0}(t_n,x) &=& \Lambda_0 t_n^{-\frac{1+\beta}{p-1}}\frac{\ln(e+(t_n^{-\frac1{2\alpha}}|x|)^2)}{1+(t_n^{-\frac1{2\alpha}}|x|)^{N+2\alpha}}
   = \Lambda_0 n^{\frac{2\alpha({1+\beta})}{p-1}}\frac{\ln(e+(n|x|)^2)}{1+(n|x|)^{N+2\alpha}} \\
  &=&\frac{\Lambda_0}{k_0} n^{\frac{2\alpha(1+\beta)}{p-1}-N}n^Nf_0(n|x|)\\
   &\ge&\frac{\Lambda_0}{k_0} {\tilde n}^{\frac{2\alpha(1+\beta)}{p-1}-N}n^Nf_0(n|x|)
  =f_{n,k_{\tilde n}}(x),
\end{eqnarray*}
where $\tilde n\le n$ and $k_{\tilde n}=\Lambda_0{\tilde n}^{\frac{2\alpha(1+\beta)}{p-1}-N}$. We see that $k_{\tilde n}=\Lambda_0{\tilde n}^{\frac{2\alpha(1+\beta)}{p-1}-N}\to\infty$ as ${\tilde n}\to\infty$, since $\frac{2\alpha(1+\beta)}{p-1}-N>0$.
Let $u_{n,k_{\tilde n}}$ be the solution of (\ref{he 2.1}) with initial data $f_{n,k_{\tilde n}}$. By Lemma \ref{lm 14-09-0},
$w_{\Lambda_0}(\cdot+t_n,\cdot)$ is a super-solution of (\ref{he 2.1}) with initial data $w_{\Lambda_0}(t_n,\cdot)$,
that is, for $(t,x)\in Q_\infty$,
$$
\partial_tw_\lambda(t+t_n,x)+(-\Delta)^\alpha w_\lambda(t+t_n,x)+(t+t_n)^{\beta}w_\lambda^p(t+t_n,x)\ge 0.
$$
By the Comparison Principle,
$$
u_{n,k_{\tilde n}}(t,x)\le w_{\Lambda_0}(t+t_n,x),\qquad \forall (t,x)\in Q_\infty,
$$
for any $\tilde n\le n$. Letting $n\to\infty$ we infer
\begin{equation}\label{18-09-0}
u_{k_{\tilde n}}(t,x)\le w_{\Lambda_0}(t,x),\qquad \forall (t,x)\in Q_\infty,
\end{equation}
where $u_{k_{\tilde n}}$ is the solution of (\ref{he 2.1}) with $k_{\tilde n}\delta_0$ initial data.
Thus (\ref{13-09-9}) is obtained by letting $\tilde n\to\infty$. Finally  (\ref{13-09-90})
 follows by the fact that
 $$\lim_{t\to0^+}w_{\Lambda_0}(t,x)=0,\qquad \forall x\in\R^N\setminus\{0\}, $$
which completes the proof.\hfill$\Box$\medskip

\medskip

\begin{lemma}\label{lm 06-10-0}
Assume $1<p<p^*_\beta$, then there exists $c_{31}>0$ such that
\begin{equation}\label{06-10-0}
u_\infty(t,x)\ge \frac{c_{31} t^{-\frac{1+\beta}{p-1}}}{1+|t^{-\frac1{2\alpha}}x|^{N+2\alpha}},\qquad \forall(t,x)\in(0,1)\times\R^N.
\end{equation}
\end{lemma}
{\bf Proof.} We divide the proof into two steps.\smallskip

\noindent {\it Step 1.} Let $\sigma_0=1+\beta-\frac{N}{2\alpha}(p-1)>0$, $\eta(t)=2-t^{\sigma_0}$ for $t>0$
 and  denote
$$v_\epsilon(t,x)=\epsilon \eta(t)\Gamma_\alpha(t,x),$$
where  $\Gamma_\alpha$ is the fundamental solution of (\ref{he 2.1}).
In this step we prove that there exists $\epsilon_0>0$ such that
\begin{equation}\label{20-10-1}
u_{k_0}\ge v_{\epsilon_0}\quad \text{ in }\; (0,1)\times\R^N,
\end{equation}
where $k_0=2\epsilon_0$ and $u_{k_0}$ is the solution of (\ref{he 2.1}) with initial data $k_0\delta_0$.
Indeed, $$\partial_t v_\epsilon(t,x)=\epsilon\eta'(t)\Gamma_\alpha(t,x)+\epsilon\eta(t)\partial_t \Gamma_\alpha(t,x)$$
and $$(-\Delta)^\alpha v_\epsilon(t,x)=\epsilon\eta(t)(-\Delta)^\alpha\Gamma_\alpha(t,x).$$
Let $\Gamma_1(t^{-\frac1{2\alpha}}x)=\Gamma_\alpha(1,t^{-\frac1{2\alpha}}x)$, then there exists $\epsilon_0>0$ such that
for any $\epsilon\le \epsilon_0$ and  $(t,x)\in(0,1)\times\R^N$, we have that
$$
\begin{array}{lll}
\displaystyle
\partial_t v_\epsilon(t,x)+ (-\Delta)^\alpha v_\epsilon(t,x)+t^\beta v_\epsilon^p(t,x)
\\[4mm]\phantom{----}\displaystyle=\epsilon\eta'(t)t^{-\frac N{2\alpha}}\Gamma_1(t^{-\frac1{2\alpha}}x)
+\epsilon^p\eta^p(t)t^{-\frac N{2\alpha}p+\beta}\Gamma_1^p(t^{-\frac1{2\alpha}}x)
\\[4mm]\phantom{----}\displaystyle
\le -\epsilon\sigma_0t^{-\frac N{2\alpha}-1+\sigma_0}\Gamma_1(t^{-\frac1{2\alpha}}x)+ 2^p\epsilon^p t^{-\frac N{2\alpha}p+\beta}\Gamma_1^p(t^{-\frac1{2\alpha}}x)
\le 0,
\end{array}
$$
the last inequality holds since $-\frac N{2\alpha}-1+\sigma_0=-\frac N{2\alpha}p+\beta$ and $\Gamma_1$ is bounded.
In particular, there holds
\begin{equation}\label{21-10-0}
\partial_t v_{\epsilon_0}(t,x)+ (-\Delta)^\alpha v_{\epsilon_0}(t,x)+t^\beta v_{\epsilon_0}^p(t,x)\le 0,\qquad \forall(t,x)\in(0,1)\times\R^N.
\end{equation}

Let  $f_n(x)=v_{\epsilon_0}(t_n,x)$ with $t_n=n^{-2\alpha}$. Since $\lim_{t\to0^+}\eta(t)=2$, then we have that
$f_n\to 2\epsilon_0\delta_0$ as $n\to\infty$ in the weak sense of measures.
There exists $N_0>0$ such that  $t_n\in (0,\frac18)$ for $n\ge N_0$. Let $w_n$ be the solution of (\ref{he 2.1}) with
initial data $f_n$, then it infers that
$$w_n(t,x)\ge v_{\epsilon_0}(t+t_n,x),\qquad (t,x)\in (0,1-t_n)\times\R^N.$$
Because  $u_{k_0}$ is uniquely defined, there holds
$$w_n\to u_{k_0}\quad \text{ as }\; n\to\infty \qquad \text{in }\; (0,1)\times\R^N$$
and
$$\lim_{n\to\infty}v_{\epsilon_0}(t+t_n,x)=v_{\epsilon_0}(t,x),\qquad \forall(t,x)\in(0,1)\times\R^N,$$
which imply (\ref{20-10-1}).\smallskip

\noindent{\it Step 2. We claim that (\ref{06-10-0}) holds.} Since
$$v_{\epsilon_0}(t,x)\ge \epsilon_0 t^{-\frac{N}{2\alpha}}\Gamma_1(t^{-\frac1{2\alpha}}x),\qquad (t,x)\in(0,1)\times\R^N,$$
 then, along with the relation $T_\lambda[u_k]=u_{k\lambda^{\frac{2\alpha(1+\beta)}{p-1}-N}}$, we observe that for any $\lambda>0$,
\begin{eqnarray*}
u_{k_0\lambda^{\frac{2\alpha(1+\beta)}{p-1}-N}}(t,x) &=& \lambda^{\frac{2\alpha(1+\beta)}{p-1}}u_{k_0}(\lambda^{2\alpha}t,\lambda x) \\
   &\ge& \lambda^{\frac{2\alpha(1+\beta)}{p-1}}v_{\epsilon_0}(\lambda^{2\alpha}t,\lambda x)
   \\&\ge &\epsilon_0\lambda^{\frac{2\alpha(1+\beta)}{p-1}-N}t^{-\frac{N}{2\alpha}}\Gamma_1(t^{-\frac1{2\alpha}}x).
\end{eqnarray*}
Let $\varrho=\lambda^{\frac{2\alpha(1+\beta)}{p-1}-N}$,
$t_\varrho=(2\varrho)^{\frac1{\frac N{2\alpha}-\frac{1+\beta}{p-1}}}$ and $T_\varrho=\varrho^{\frac1{\frac N{2\alpha}-\frac{1+\beta}{p-1}}}$, then
$$0<t_\varrho<T_\varrho\to0\quad{\rm as}\ \ \varrho\to\infty.$$
For $(t,x)\in(t_\varrho,T_\varrho)\times\R^N$, we have that
\begin{eqnarray*}
 u_{k_0\varrho}(t,x) \ge \epsilon_0\varrho t^{-\frac{N}{2\alpha}}\Gamma_1(t^{-\frac1{2\alpha}}x) \ge \frac{\epsilon_0}{2} t^{-\frac{1+\beta}{p-1}}\Gamma_1(t^{-\frac1{2\alpha}}x),
\end{eqnarray*}
then
$$u_\infty(t,x) \ge  \frac{\epsilon_0}{2} t^{-\frac{1+\beta}{p-1}}\Gamma_1(t^{-\frac1{2\alpha}}x),\qquad \forall (t,x)\in(t_\varrho,T_\varrho)\times\R^N.$$
which implies (\ref{06-10-0}) and completes the proof. \hfill$\Box$

\medskip

\noindent{\bf Proof of Theorem \ref{teo 3}.} It follows from Proposition \ref{lm 13-09-1} and Lemma \ref{lm 13-09-2} that
$u_\infty$ is a nontrivial self-similar solution of (\ref{he 2.1}) and (\ref{13-09-08}) follows by (\ref{13-09-9}), (\ref{06-10-0}) and
$\ln(e+|t^{-\frac1{2\alpha}}x|^2)\le 2\ln(2+|t^{-\frac1{2\alpha}}x|)$, which ends the proof. \hfill$\Box$\medskip

We have actually a stronger result which is a consequence of Theorem \ref{teo 5}-(i) proved in next section:

\begin{corollary}\label{coro} Assume $1+\frac{2\alpha(1+\beta)}{N+2\alpha}<p<1+\frac{2\alpha(1+\beta)}{N}$. Then \smallskip

\noindent either
\begin{equation}\label{6.1}
  \tilde u> u_\infty\quad {\rm in}\quad Q_\infty
\end{equation}
or
\begin{equation}\label{6.2}
\tilde u\equiv u_\infty\quad {\rm in}\quad Q_\infty.
\end{equation}

\end{corollary}
\subsection{The case $1<p<1+\frac{2\alpha(1+\beta)}{N+2\alpha}$}

For $1<p<1+\frac{2\alpha(1+\beta)}{N+2\alpha}$, it follows from Lemma \ref{lm 06-10-0} that
\begin{equation}\label{08-10-0}
  \lim_{t\to0^+}u_\infty(t,x)=\infty,\qquad \forall x\in\R^N.
\end{equation}

\noindent{\bf Proof of Theorem \ref{teo 4} $(i)$.}
Let $f_0\in C_c(\R^N)$ be a nonnegative function such that
$${\rm supp}f_0\subset B_1(0)\quad \ {\rm and} \quad \ \max_{x\in B_1(0)}f_0=1.$$
Denote
$$f_{n,k}(x)=kn^{\theta N}f_0(n^\theta(x-x_0)),$$
where $k\le n^\tau$ with $\tau=\frac12(\frac{2\alpha(1+\beta)}{p-1}-N-2\alpha)>0$, $\theta=\frac{\tau}{N}$ and $x_0\in\R^N$.
Since $f_{n,k}(x)\le n^{\tau}$ for $x\in B_1(x_0)$, $f_n(x)=0$ for $x\in B_1^c(x_0)$ and
$$
v_{\epsilon_0}(t_n,x) \ge\frac{c_{39} n^{\frac{2\alpha(1+\beta)}{p-1}-N-2\alpha}}{(2+|x_0|)^{N+2\alpha}},\qquad \forall x\in B_1(x_0),
$$
where $t_n=n^{-2\alpha}$. Then there exists $N_0>0$ such that for any $n\ge N_0$,
 $$f_{n,k}(x)\le v_{\epsilon_0}(t_n,x),\qquad \forall x\in B_1(x_0).$$
Since $n^{\theta N}f_0(n^\theta(x-x_0))\to c_{41}\delta_{x_0}$,  as $ n\to\infty$ in weak sense of measures, for some $c_{41}>0$.
\smallskip

Let $w_{n,k}$ be the solution of (\ref{he 2.1}) with initial data $f_{n,k}$, then
$$w_{n,k}(0,x)=f_{n,k}(x)\le v_{\epsilon_0}(t_n,x)\le u_\infty(t_n,x),\qquad \forall x\in\R^N.$$
Therefore, by the Comparison Principle $$w_{n,k}(t,x)\le  u_\infty(t+t_n,x),\qquad \forall (t,x)\in Q_\infty.$$
We observe that
$$\lim_{k\to\infty}[\lim_{n\to\infty}w_{n,k}(t,x)]= u_\infty(t,x-x_0),\qquad \forall (t,x)\in Q_\infty.$$
Thus, we derive that
\begin{equation}\label{08-10-1}
 u_\infty(t,x-x_0)\le u_\infty(t,x),\qquad \forall (t,x)\in Q_\infty.
\end{equation}
Then
$ u_\infty(t,x-x_0)= u_\infty(t,x)$ for all $(t,x)\in Q_\infty,$
which  implies that $u_\infty$ is independent of $x$.
Combining  (\ref{02-10-3}) and (\ref{06-10-0}), implies that
$$u_\infty=\left(\frac{1+\beta}{p-1}\right)^{\frac{1}{p-1}}t^{-\frac{1+\beta}{p-1}}.$$
The proof is complete.\hfill$\Box$\\[2mm]

In the case of $p=1+\frac{2\alpha(1+\beta)}{N+2\alpha}$, it derive from Lemma \ref{lm 06-10-0} that
$$
\liminf_{t\to0^+}u_\infty(t,x)\geq
  \lim_{t\to0^+}\frac{c_{40} t^{-\frac{1+\beta}{p-1}}}{1+|t^{-\frac1{2\alpha}}x|^{N+2\alpha}}
  =\frac{c_{40}}{|x|^{N+2\alpha}},\qquad \forall x\in\R^N.
$$

\noindent{\bf Proof of Theorem \ref{teo 4} $(ii)$.} We note that
$u_\infty$ is a self-similar solution  of (\ref{he 2.1}). Moreover, we derive (\ref{08-10-2})  by (\ref{06-10-0}), which
ends the proof.   \hfill$\Box$

\subsection{The self-similar equation}
In this section we prove Theorem \ref{teo 5}.\medskip

\noindent{\bf Proof of Theorem \ref{teo 5} $(i)$.} We set $v_\infty(\eta)=t^{\frac{1+\beta}{p-1}}u_\infty(1,\eta)$. Then relations (\ref{decay}) and (\ref{decay2})
hold from Lemmas \ref{lm 13-09-2} and \ref{lm 06-10-0}. Assume $\tilde v$ is another positive solution of (\ref{13-09-07}). Then
$(t,x)\mapsto t^{-\frac{1+\beta}{p-1}}\tilde v(t^{-\frac{1}{2\alpha}}x)$ is a positive self-similar solution of (\ref{5.1}). By Proposition \ref{lm 5.1} it is larger than $u_\infty$. Thus $v_\infty\leq \tilde v$. Assume now that there exists $\eta_0\in \R^N$ such that
$v_\infty(\eta_0)=\tilde v(\eta_0)$. and set $w=\tilde v-v_\infty$. Then
$$\begin{array} {lll}\displaystyle
(-\Delta)^\alpha w(\eta_0)=\lim_{\epsilon\to 0}(-\Delta)_\epsilon^\alpha w(\eta_0)
=\lim_{\epsilon\to 0}\myint{B^c_\epsilon(\eta_0)}{}\myfrac{w(\eta_0)-w(\eta)}{|\eta-\eta_0|^{N+2\alpha}}d\eta
<0.
\end{array}$$
Since $\nabla w(\eta_0)$ we reach a contradiction.    \hfill$\Box$

\medskip

\noindent{\bf Proof of Theorem \ref{teo 5} $(ii)$.} It is a consequence of the equality
$$u_\infty=U_p\Longleftrightarrow v_\infty=\left(\frac{1+\beta}{p-1}\right)^{\frac{1}{p-1}}.$$
\medskip

\noindent{\bf Open problem.} We conjecture that in the case $1+\frac{2\alpha (1+\beta)}{N+2\alpha}<p<1+\frac{2\alpha (1+\beta)}{N}$,
$v_\infty$ is the unique positive solution of the self-similar equation satisfying (\ref{decay}). One step could be to prove that any positive solution $\tilde v$ satisfying (\ref{decay}) satisfies, for some $K>1$,

\begin{equation}\label{asympt0}
\tilde v\leq K v_\infty\qquad\text{in }\;\R^N.
\end{equation}
We also conjecture that $v_\infty$ satisfies the following asymptotic behavior
\begin{equation}\label{asympt}
 v_\infty(\eta)=c_{N,p,\alpha,\beta}|\eta|^{-N-2\alpha}\qquad\text{as }\;|\eta|\to \infty.
\end{equation}
Thus if any  positive solution $\tilde v$ inherits the same property, the conclusion (and the uniqueness)   follows.



\bigskip

\noindent{\bf Acknowledge:} H. Chen is supported by National Natural Science Foundation
of China,  No:11401270 and the Project-sponsored by SRF for ROCS, SEM.
L. V\'{e}ron is supported by  the MathAmsud collaboration
program 13MATH-02 QUESP.

\end{document}